\documentclass{article}

    \PassOptionsToPackage{numbers, compress}{natbib}

    \usepackage[preprint]{neurips_2020}

\usepackage{hyperref}

\usepackage[utf8]{inputenc} 
\usepackage[T1]{fontenc}    
\usepackage{booktabs}       
\usepackage{amsfonts}       
\usepackage{nicefrac}       
\usepackage{microtype}      
\usepackage{bm}
\usepackage{mathtools}
\usepackage{mathrsfs}
\usepackage{graphicx}
\usepackage{algorithm}
\usepackage{algpseudocode}
\usepackage[nolist]{acronym}
\usepackage{bbm}
\usepackage{xcolor}
\usepackage{wrapfig}
\usepackage{siunitx}
\usepackage{multirow}
\usepackage{enumitem}

\usepackage[symbol]{footmisc}
\usepackage{float}
\usepackage{cleveref}

\usepackage{natbib}

\newcommand{\rd}{{\mathrm d}}

\newcommand{\calU}{{\cal{U}}}

\renewcommand{\labelenumi}{\roman{enumi})}

\newcommand{\calN}{{\cal N}}

\newcommand{\argmax}{\operatornamewithlimits{argmax}}
\newcommand{\argmin}{\operatornamewithlimits{argmin}}
\newcommand{\T}{^\mathsf{T}}

\newcommand{\Rb}{\mathbb{R}}
\newcommand{\Eb}{\mathbb{E}}
\newcommand{\CVaR}{\text{CVaR}}
\newcommand{\VaR}{\text{VaR}}

\newcommand{\val}[2]{\mathit{V}_{#2}\left({#1}\right)}

\begin{acronym}
\acro{CVaR}{Conditional Value-at-Risk}
\acro{PMP}{Pontryagin Maximum Principle}
\acro{VaR}{Value-at-Risk}
\acro{GASS}{Gradient-based Adaptive Stochastic Search}
\acro{SA}{Stochastic Approximation}
\acro{SOC}{Stochastic Optimal Control}
\acro{RL}{Reinforcement Learning}
\acro{SDPG}{Sample-based Distributional Policy Gradients}
\acro{RS3}{Risk Sensitive Stochastic Search}
\acro{CEM}{Cross Entropy Method}
\acro{MPC}{Model Predictive Control}
\acro{SS}{Stochastic Search}
\acro{MPPI}{Model Predictive Path Integral}
\end{acronym}

\DeclarePairedDelimiter\abs{\lvert}{\rvert}%
\DeclarePairedDelimiter\norm{\lVert}{\rVert}%

\makeatletter
\let\oldabs\abs
\def\abs{\@ifstar{\oldabs}{\oldabs*}}
\let\oldnorm\norm
\def\norm{\@ifstar{\oldnorm}{\oldnorm*}}
\makeatother

\renewcommand{\thefootnote}{\fnsymbol{footnote}}

\title{Adaptive Risk Sensitive Model Predictive Control with Stochastic Search}

\author{Ziyi Wang\thanks{Correspondence to: \texttt{ZiyiWang@gatech.edu}},~ Oswin So, ~ Keuntaek Lee, ~ Camilo A. Duarte, ~ Evangelos  A. Theodorou \\
Autonomous Control and Decision Systems Laboratory \\
School of Aerospace Engineering\\
Georgia Institute of Technology
}

\begin{document}

\maketitle
\begin{abstract}%
We present a general framework for optimizing the Conditional Value-at-Risk for dynamical systems using path space stochastic search. The framework is capable of handling the uncertainty from the initial condition, stochastic dynamics, and uncertain parameters in the model. The algorithm is compared against a risk-sensitive distributional reinforcement learning framework and demonstrates improved performance on a pendulum and cartpole with stochastic dynamics. We also showcase the applicability of the framework to robotics as an adaptive risk-sensitive controller by optimizing with respect to the fully nonlinear belief provided by a particle filter on a pendulum, cartpole, and quadcopter in simulation.
\end{abstract}

\section{Introduction}
\begin{wrapfigure}{r}{0.65\textwidth}
    \vspace{-10 pt}
    \centering
    \includegraphics[width=0.65\textwidth]{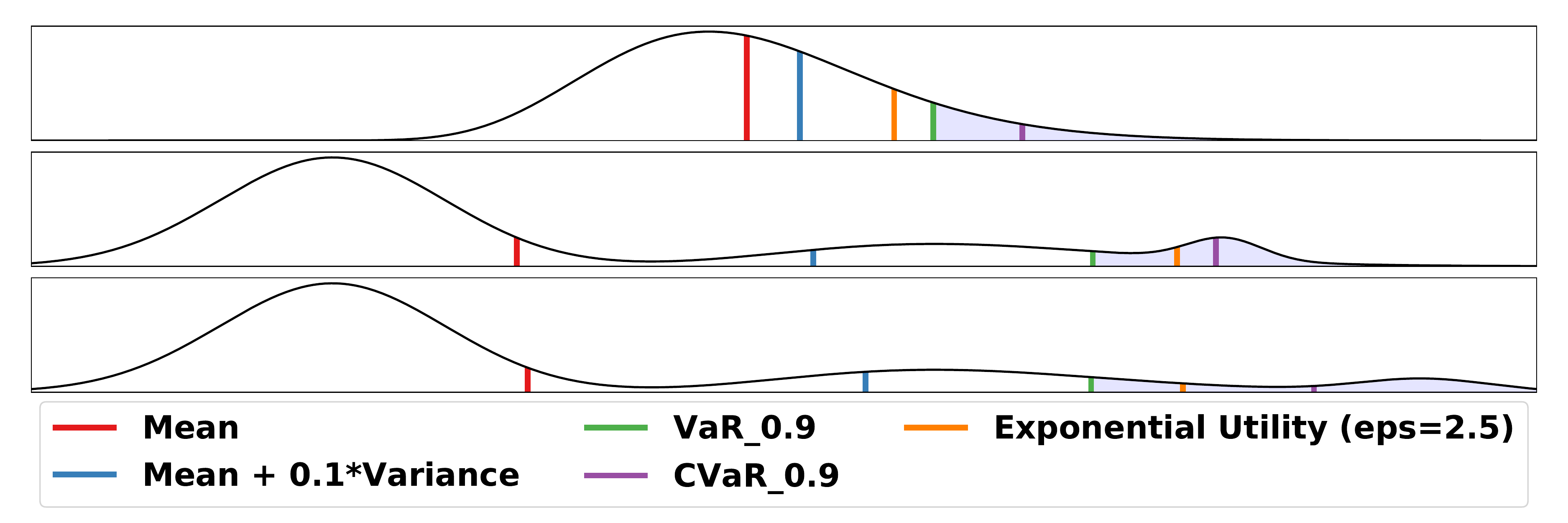}
    \vspace{-20 pt}
    \caption{
    Comparison of different risk measures across three different distributions of costs.
    }
    \label{fig:cvar_graph}
    \vspace{-10 pt}
\end{wrapfigure}

The majority of \ac{SOC} and \ac{RL} literature handles uncertainty in dynamical optimization problems by simply optimizing with respect to the expected cost/reward. In many applications, however, it is desirable to consider the risk associated with a policy instead of its performance on average. A simple but practical risk measure is the variance or mean-plus-variance \cite{gosavi2014variance,di2012policy}. A major problem with variance is that it is a symmetric risk measure. The undesired high cost scenario is penalized the same way as the desired low cost outcome. Common asymmetric risk measures include exponential utility \cite{howard1972risk, fleming1995risk}, which quantifies the exponential growth of risk as the cost increases, and \ac{VaR} \cite{jorion2000value}, $\text{VaR}^\gamma(X) = \inf\{t:\mathbb{P}(X\leq t)\geq \gamma\}$, which quantifies statistically the $\gamma$-quantile of the uncertain cost distribution with $\gamma\in(0,1)$ being the risk level. While exponential utility and \ac{VaR} penalize one side of the cost distribution as desired, they are not \textit{coherent} risk measures (see 
Appendix \ref{sec:coherency}
for definition) \cite{artzner1999coherent}. \ac{CVaR} \cite{rockafellar2000optimization} is a natural extension to \ac{VaR} defined as \begin{equation*}
\CVaR^\gamma(X) = \frac{1}{1-\gamma}\int_\gamma^1 \text{VaR}^r(X)\mathrm{d}r
\end{equation*}
which is equivalent to the conditional expectation beyond \ac{VaR}, $\mathbb{E}[X|X\geq \text{VaR}^\gamma(X)]$, if $X$ has continuous distribution. The main advantages of \ac{CVaR} are that it is coherent, measures only the worst cases compared to exponential utility, but takes into account the entire tail instead of only the $\gamma$-quantile compared to \ac{VaR}. Figure \ref{fig:cvar_graph} illustrates the difference between risk measures on three distributions. Comparing the top and middle distributions, it is clear that symmetric risk measures cannot capture the risk associated with a heavy tail. From the middle and bottom distributions, it can be observed that \ac{CVaR} is more sensitive to the tail distribution than \ac{VaR}.

\ac{CVaR} has been used as a risk metric extensively in the field of finance \cite{agarwal2004risks, krokhmal2002portfolio,zhu2009worst}, power utility \cite{conejo2010decision, morales2010short}, supply chain management \cite{chen2007risk}, etc. In recent years, it is also seeing a rise in popularity in robotics research. However, \ac{CVaR} optimization for dynamical systems suffers from the problem of time inconsistency \cite{bjork2010general}, meaning that the optimal policy at a particular timestep might be suboptimal at a future time. We encourage the readers to refer to \cite{shapiro2014lectures} Section 6.8.5 for the mathematical definition of time consistency and \cite{chow2015time} Example 2 for an intuitive example on the time inconsistency of \ac{CVaR}. The time inconsistency makes directly applying popular methods in \ac{SOC} and \ac{RL} that originated from dynamic programming to the \ac{CVaR} optimization problem hard. Several methods have been proposed to overcome this problem by lifting the state space of the problem. In \cite{pflug2016time} and \cite{chow2015risk}, the \ac{CVaR} decomposition theorem is introduced to obtain a dual representation of \ac{CVaR} and optimizes the associated Bellman equation over a space of probability densities. Alternatively, the convex extremal formulation of \ac{CVaR} \cite{rockafellar2000optimization} can be used to alleviate the time-inconsistency problem \cite{bauerle2011markov}. Both approaches, however, require some form of state space augmentation and require solving an additional optimization problem. On the other hand, algorithms that directly optimize the policy \cite{tamar2014policy} are not affected by the time-inconsistency problem since they do not rely on the dynamic programming principles.

A promising sampling-based approach that directly optimizes the policy for solving general nonlinear optimization problems is stochastic search. Stochastic search is a general class of optimization methods that optimizes the objective function by randomly sampling and updating candidate solutions. Many well-known algorithms, such as \ac{CEM} \cite{rubinstein2013cross}, genetic algorithm \cite{zames1981genetic}, and simulated annealing \cite{kirkpatrick1983optimization} fall into this category. Recently, a \ac{GASS} algorithm was proposed by \cite{zhou2014gradient}. \ac{GASS} updates the candidate solution by taking the gradient with respect to the sampling distribution parameters of solutions and approximating the gradient with Monte Carlo sampling. \cite{boutselis2019constrained} extended \ac{GASS} to constrained dynamic optimization problems, and \cite{wang2019information} showed that Information Theoretic \ac{MPPI} \cite{williams2017information} emerges as a special case of dynamic \ac{GASS} in the case of Gaussian and Poisson sampling distributions.

In this paper, we extend the risk sensitive formulation of \ac{GASS} \cite{zhu2018simulation} and present \ac{RS3}, a general framework for solving \ac{CVaR} optimization for dynamical systems. The resulting algorithm bypasses the problem of time-inconsistency by directly performing stochastic gradient descent on the sampling distribution parameters. The framework is capable of handling uncertain initial states, parameters, and system stochasticity. We demonstrate the applicability of our framework to robotics by combining it with a particle filter to perform risk-sensitive belief space control on a pendulum, cartpole and quadcopter in simulation. In addition, we show that our framework outperforms in terms of final cost in average and under other risk measures compared to a risk-sensitive distributional \ac{RL} algorithm, the \ac{SDPG} \cite{singh2020cvarsdpg} on a pendulum and cartpole system in simulation.

The rest of this paper is organized as follows: in \Cref{sec:problem_formulation} we formulate the problem and derive the stochastic search framework. We present the algorithm and its application to belief space optimization in \Cref{sec:algorithm}. The simulation results are included in \Cref{sec:results}. Finally, we conclude the paper in \Cref{sec:conclusion}.

\section{Stochastic Search}
\label{sec:problem_formulation}

\subsection{Notation, Mathematical Preliminaries, and Problem Formulation}
Note that hereon we use $\mathbb{E}_{p(x)}[f(x)]$ to denote the integral $\int_{\Omega_x} f(x) p(x) \mathrm{d}x$ and $p(x)$ is dropped from the expectation for simplicity when it is clear which distribution the expectation is taken with respect to. We consider the problem of minimizing the \ac{CVaR} of cost function $J:\mathbb{R}^{{n_x}\times T}\times \calU \to \mathbb{R}^+$
\begin{equation}
    U^* = \argmin_{U\in \calU}\CVaR^{\gamma}[J(X,U)],
    \label{eq:cvar_opt}
\end{equation}
subject to nonlinear stochastic dynamics
\begin{equation}
    x_{t+1} \sim p(x_{t+1}|x_t,u_t; \phi).
    \label{eq:dynamics}
\end{equation}
Here we have $U=\{u_0,\cdots,u_{T-1}\}\in\calU$ as the control path and $X=\{x_0,\cdots,x_T\}\in\Rb^{n_x\times T}$ as the state path where $T\in[0,\infty)$ is the optimization horizon. $\calU\subset\Rb^{n_u\times T-1}$ is the set of admissible control sequences, and $\phi$ is the system parameters. This formulation is capable of handling stochasticity in dynamics, $p(x_{t+1}|x_t,u_t)$, uncertain parameters, $p(\phi)$, and uncertain initial condition, $p(x_0)$. The \ac{CVaR} is computed with respect to the uncertainty distributions. Assuming that $p(x_{t+1}|x_t,u_t;\phi)$, $p(x_0)$ and $p(\phi)$ are independent continuous density functions and $J$ is continuous, the minimization problem \eqref{eq:cvar_opt} can be rewritten as
\begin{equation}\label{eq: before_sample_dist}
    U^* = \argmin_{U\in \calU}\Eb_{p(\chi)}[J(X,U)|J\geq \VaR^\gamma(J)].
\end{equation}
where $p(\chi)=p(x_0)p(\phi)\prod_{t=0}^T p(x_{t+1}|x_t,u_t)$ is the joint pdf of all uncertainty distributions. We parameterize the control $u$ with a policy $\pi_\eta$ characterized by its parameters $\eta$. The policy can be of any functional form, i.e. open-loop ($u_t = v_t, \eta_t = v_t$), linear feedback ($u_t = k_tx_t + v_t, \eta_t = \{k_t, v_t\}$) or a neural network. We then define a sampling distribution for the policy parameters from the exponential family with a pdf of the form
\begin{equation}
p(\eta_t;\theta_t) = h(\eta_t)\exp\Big(\theta_t^\mathrm{T}T(\eta_t)-A(\theta_t)\Big),
\end{equation}
where $\theta$ is the natural parameters of the distribution and $T(\eta)$ is the sufficient statistics of $\eta$. The minimization is now performed with respect to the natural parameters
\begin{equation} \label{eq:after_sample_dist}
    \theta^* = \argmin_{\theta\in \Theta}\Eb_{p(\chi,\theta)}[J(X,U)|J\geq \VaR^\gamma(J)].
\end{equation}
The expectation is now taken with respect to the joint distribution of uncertainty and the sampling distribution. It is easy to show that the expected cost in \eqref{eq:after_sample_dist} is an upper bound on the optimal cost in \eqref{eq: before_sample_dist}.

\subsection{Update Law Derivation}
In this section we derive the gradient descent update rule for the optimization problem defined in \eqref{eq:after_sample_dist}. To be consistent with \ac{GASS}, we turn the minimization problem into a maximization one by optimizing with respect to $-J$. We then introduce a shape function $S:\Rb \rightarrow \Rb^+$ which allows for different weighing schemes of the cost levels, leading to different optimization behaviors \cite{ollivier2017information}. The problem is then transformed into
\begin{equation}\label{eq:gass_obj}\begin{split}
    \theta^* &= \argmax_{\theta\in \Theta}\Eb_{p(\theta)}[S(-\Eb_{p(\chi)}[J(X,U)|J\geq \VaR^\gamma(J)])]\\
    &=\argmax_{\theta\in \Theta}\Eb_{p(\theta)}[S(-\CVaR^\gamma[J(X,U)])].
    \end{split}
\end{equation}
The shape function needs to satisfy the following conditions:
\begin{enumerate}
    \item $S(y)$ is nondecreasing in $y$ and bounded from above and below for bounded $y$, with the lower bound being away from zero.
    \item The set of optimal solutions $\{\argmax_{y\in\mathcal{Y}} S(H(y))\}$ after the transform is a non-empty subset of the solutions $\{\argmax_{y\in\mathcal{Y}} H(y)\}$ of the original problem.
\end{enumerate}
Common shape functions include: 1. the identity function, $S(y)=y$; 2. the exponential function, $S(y;\kappa)=\exp(\kappa y)$, which leads to the Information Theoretic \ac{MPPI} update law \cite{williams2017information}; 3. the sigmoid function, $S(y; \kappa, \varphi)=(y-y_{\text{lb}})\tfrac{1}{1+\exp(-\kappa(y-\varphi))}$, where $y_{\text{lb}}$ is a lower bound for the cost and $\varphi$ is the ($1-\rho$)-quantile, which results in an updste law similar to \ac{CEM} but with soft elite threshold $\rho$.

Finally, we apply another log transformation to obtain scale-free gradient and the optimization problem becomes
\begin{equation}\label{eq:gass_log_obj}
    \theta^* = \argmax_{\theta\in \Theta}\ln\Eb[S(-\CVaR^\gamma[J(X,U)])]=\argmax_{\theta\in \Theta} l(\theta).
\end{equation}
Since $\ln:\Rb^+\rightarrow \Rb$ is a strictly increasing function, it does not change the maximization objective. We can now take its gradient with respect to the parameters. Writing the expectation as an integral with respect to the path probability $p(X,U;\theta) $  we get
\begin{equation}\label{eq:cvar_integral}
    \Eb[S(-\CVaR^\gamma[J(X,U)])]=\int_{\Omega_\eta} S\Big(-\int_{\Omega_\chi}
    J(X,U)p(X,U;\theta)\rd \chi\Big)\rd \eta,
\end{equation}
where $\Omega_\chi$ is defined such that $J\geq \VaR^{\gamma}(J)$ if and only if $\chi\in\Omega_\chi$, and $\Omega_\eta$ is defined such that $\pi_{\eta_t}(x_t)\in\calU_t, \forall t$. The path probability distribution can be decomposed as
\begin{align}
    p(X,U;\theta) &= p(x_T|x_{T-1},\pi_\eta(x_{T-1});\phi)p(\eta_{T-1};\theta_{T-1})\cdots p(x_0)p(\phi)\\
    &= \underbrace{p(x_0)p(\phi)\prod_{t=0}^{T-1}p(x_{t+1}|x_t,\pi_\eta(x_t);\phi)}_{p(\chi)}\underbrace{\prod_{t=0}^{T-1}p(\eta_t;\theta_t)}_{p(\eta)}.
\end{align}
Note that since the uncertainty and sampling distribution are independent, their joint distribution can be broken into the product of the two. The gradient of the objective function \eqref{eq:gass_log_obj} with respect to the parameters can be taken as
\begin{equation}\label{eq:gass_update_law}
    \nabla_{\theta} l(\theta) = \frac{\Eb[S(-\CVaR^\gamma[J(X,U)])\nabla_{\theta}\big(\sum_{t=0}^{T-1}\ln p(\eta_t;\theta_t)\big)]}{\Eb[S(-\CVaR^\gamma [J(X,U)])]}.
\end{equation}
\Cref{sec:derivation_appendix}
. The gradient of the log parameter distribution at each time step can be calculated as
\begin{align}
    \nabla_{\theta_t} \ln p(\eta_t;\theta_t) &= \nabla_{\theta_t} \ln \Big( h(\eta_t)\exp\Big(\theta_t^\mathrm{T}T(\eta_t)-A(\theta_t)\Big)\Big) \\
    &= \nabla_{\theta_t} (\theta_t^\mathrm{T}T(\eta_t)-A(\theta_t))\\
    &= T(\eta_t) - \nabla_{\theta_t} A(\theta_t).
\end{align}
Plugging it back into the gradient of the cost function, we get
\begin{equation}
    \nabla_{\theta_t} l(\theta) = \frac{\Eb[S(-\CVaR^\gamma[J(X,U)])(T(\eta_t)-\nabla_{\theta_t}A(\theta_t))]}{\Eb[S(-\CVaR^\gamma [J(X,U)])]}.
    \label{eq:update_law1}
\end{equation}
With this, we have a gradient ascent update law for the parameters as
\begin{equation}
    \theta_t^{k+1} = \theta_t^k + \alpha^k \nabla_{\theta_t} l(\theta^k).
    \label{eq:update_law2}
\end{equation}
where the step size sequence $\alpha^k$ satisfies the typical assumptions in \ac{SA}:
\begin{align}
    \alpha^k > 0 \quad \forall k, \quad
    \lim_{k \to \infty} \alpha^k = 0, \quad
    \sum_{k=0}^\infty \alpha_k = \infty
\end{align}

\subsection{Practical Considerations}
\textbf{Numerical Approximation:} The \ac{CVaR} of the cost can be approximated \cite{kolla2019concentration} with
\begin{align}
    \hat{C}_n^\gamma &= \hat{V}_n^\gamma + \frac{1}{M(1-\gamma)}\sum_{m=1}^M J_{n,m} - \hat{V}_n^\gamma\big)^+\label{eq:cvar_approx} \\
    \hat{V}_n^\gamma &= \inf \{x:\frac{1}{M}\sum_{m=1}^M\mathbbm{1}_{\{J_{n,m}\leq x\}}\geq \gamma\}.\label{eq:var_approx}
\end{align}
The outer expectation can be approximated as $\Eb[S(-\CVaR^\gamma[J(X,U)])]=\frac{1}{N}\sum_{n=1}^N S(-\hat{C}^\gamma_n)$. Note that the expectation and \ac{CVaR} in \eqref{eq:update_law1} are computed as averages over costs defined on entire trajectory samples.

\noindent \textbf{Model Predictive Control Formulation:} The parameter update in  \Cref{eq:update_law2} can be used for trajectory optimization as well as in a receding horizon or \ac{MPC} fashion. \ac{MPC} is a powerful algorithmic approach to nonlinear feedback control which is essential in tasks that involve risk measures or high order statistical characteristics of cost functions. In this paper we will leverage parallelization using GPUs to implement \Cref{eq:update_law2} in \ac{MPC} fashion.   

\noindent \textbf{Adaptive Stochastic Search:} The \ac{MPC} formulation allows online interaction with the stochastic system dynamics. Data from this online interaction  can be used to feed  adaptive  or state estimation schemes that update the probability distribution $p(\chi)$ in an online fashion. In this work, we make use of a nonlinear state estimator, namely a particle filter, to  propagate and update distribution $p(\chi)$ over time. The resulting control architecture is a sampling-based risk-sensitive adaptive MPC scheme that optimizes \ac{CVaR} while adapting the probability distribution over parametric uncertainties. The details of this approach are further explained in the next section. 

\section{Algorithm}
\label{sec:algorithm}

In this section we present the \ac{RS3} algorithm implemented in \ac{MPC} fashion, as shown in Algorithm 1. 
At initial time, the policy parameter distribution's natural parameters are initialized. Given an initial state and parameter distribution provided by an estimator, $M$ i.i.d. samples are obtained.
$N$ policies are sampled from the policy parameter distribution and each copy of the policy is applied to all $M$ samples of the initial states. In the case of stochastic dynamics, the states of each of the $M$ samples are propagated with an independent realization of the stochastic dynamics. A cost is then calculated for each of the total $N\times M$ trajectories. For each policy sample, its associated \ac{CVaR} cost is approximated with the $M$ cost samples using \eqref{eq:cvar_approx} and \eqref{eq:var_approx}. Using the \ac{CVaR} values, the policy parameter distributions' natural parameters can be updated using \eqref{eq:update_law1} and \eqref{eq:update_law2}. In our simulations, we use Gaussian distributions with fixed variance to sample policy parameters, for which the sufficient statistics are $T(\eta_t)=\eta_t$ and $\nabla_{\theta_t} A(\theta_t)=\Eb[T(\eta_t)]$. The parameter update step is detailed in algorithm 2. As is common in \ac{SA} algorithms, Polyak averaging is performed on the natural parameters to improve the convergence rate \cite{polyak1992acceleration}. With the Polyak averaged natural parameters $\bar{\eta}$, an optimal policy can be sampled and applied to the system for $\tau$ timesteps. Finally, we apply a shift operator $\omega(\theta, \tau)$ that recedes the optimization horizon and outputs $\tilde{\theta}_t = \theta_{t+\tau}$. The last $\tau$ timesteps of the natural parameters are re-initialized.


\begin{figure}[t]
    \centering
    \includegraphics[width=\textwidth]{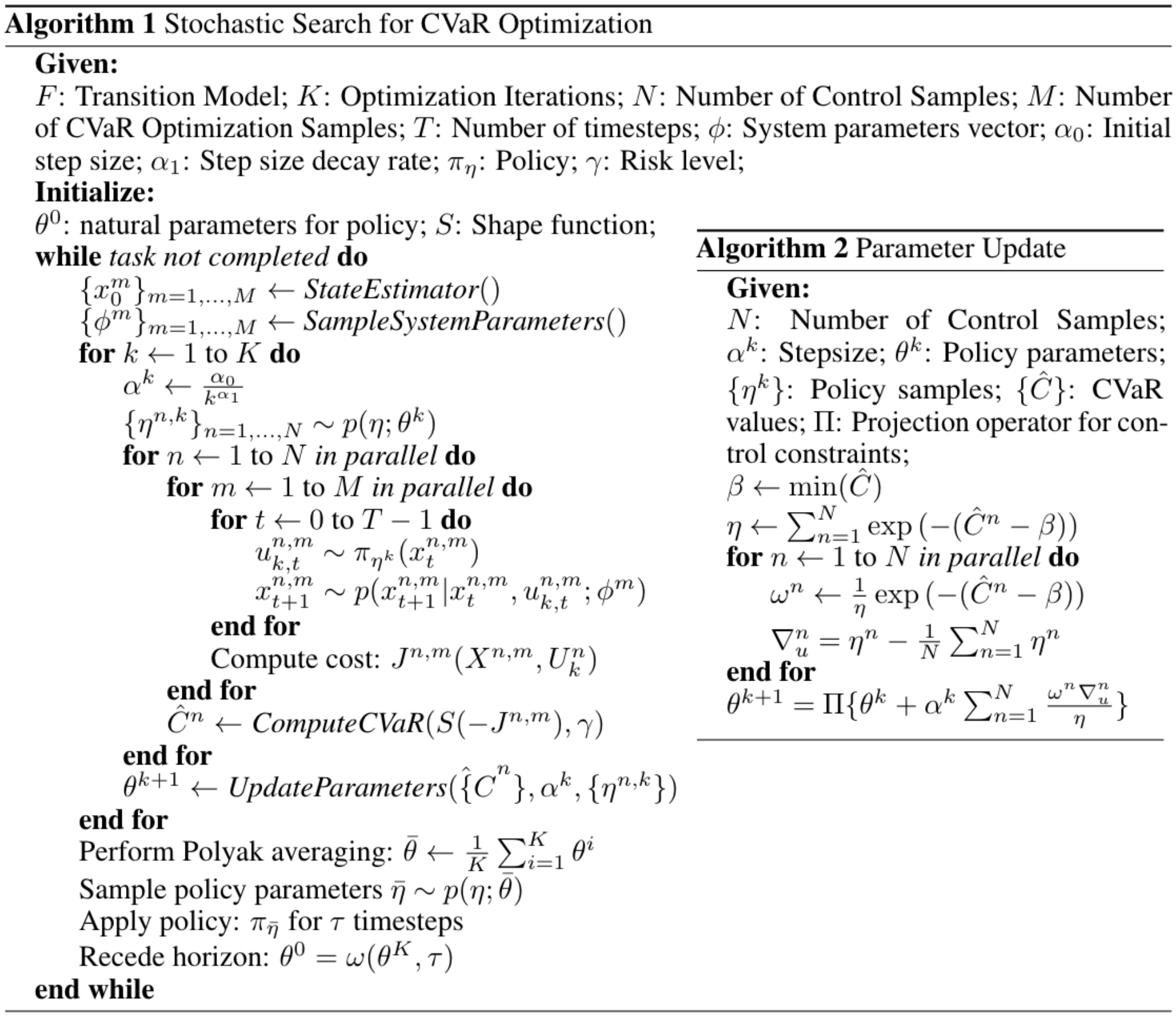}
\vspace{-0.6cm}
\label{fig:algos}
\end{figure}

\setlength{\intextsep}{10pt} 
\setlength{\floatsep}{10pt} 
\setlength{\textfloatsep}{5pt} 

Note that the \ac{RS3} algorithm can handle any or all uncertainties from initial state distribution, uncertain parameters and stochastic dynamics provided that i.i.d. samples can be generated from the uncertain distribution.

In our simulation examples of belief space control, we use a particle filter to provide the initial state distribution. To handle uncertain model parameters, we augment the states by the uncertain parameters and use a particle filter to learn its distribution (detailed in 
\Cref{sec:pf_implementation})
. In both cases, i.i.d. particles from non-Gaussian distributions can be used directly used in the \ac{RS3} algorithm. However, we want to stress that any filter can be used together with the \ac{RS3} algorithm.

\section{Simulation Results}
\label{sec:results}

In this section, we showcase the general applicability of \ac{RS3} in dealing various types of uncertainty.
1) \textit{External noise:} Comparing its performance against the \ac{SDPG} algorithm for \ac{CVaR} optimization.
2) \textit{Uncertain system parameters:} combining it with a particle filter to perform risk sensitive control in belief space. 
3) \textit{Uncertain initial condition}: This can be found in \Cref{sec:appendix_comparison_sdpg}.

All simulations were performed with a risk level of $0.9$.
The open loop policy $\pi_\eta(x) = \eta$ is used for all simulations, where $\eta$ directly maps to the controls.
The multivariate normal distribution is chosen as the sampling distribution, and we use the method proposed in \cite{boutselis2019constrained} to handle the box control constraints in the simulation tasks by sampling from a truncated multivariate normal distribution.
The tuning parameters of all simulations are included in the Appendix.

\subsection{Comparison Against \acl{SDPG}}
The \acl{SDPG} algorithm \cite{singh2020sdpg, singh2020cvarsdpg} is one of the most recent work on optimizing \ac{CVaR} for dynamical systems.
\ac{SDPG} \cite{singh2020sdpg} and the risk-sensitive version of \ac{SDPG} \cite{singh2020cvarsdpg} are actor-critic type policy gradient algorithms in the distributional \ac{RL} \cite{bellemare2017dirl} setting.

The actor network parameterizes the policy and the critic network learns the return distribution by reparameterizing simple Gaussian noise samples.
The risk-sensitive version of \ac{SDPG} \cite{singh2020cvarsdpg} is an extension of the naive \ac{SDPG} \cite{singh2020sdpg} algorithm by using \ac{CVaR} as a loss function to train the actor network to learn a risk-sensitive policy.
We compared \ac{RS3} and the risk-sensitive \ac{SDPG} on two classic control systems in OpenAI Gym \cite{brockman2016openai}, a pendulum and a cartpole.

\begin{table}[t]
  \begin{center}
    \caption{Reward comparison results between \ac{RS3} (ours) vs. risk-sensitive \ac{SDPG} \cite{singh2020cvarsdpg}}
    \label{table:comparison_sdpg}
    \begin{tabular}{|c|c|c|c|c|c|c|c|c|}
      \toprule
      \multicolumn{2}{|c|}{System} & \multicolumn{4}{c|}{Pendulum} & \multicolumn{3}{c|}{Cartpole} \\
      \hline
      \multicolumn{2}{|c|}{Noise variance} & 0.3 & 1.0 & 2.0 & 3.0 & 0.3 & 1.0 & 2.0 \\
      \hline
      & Mean                   & 169.2 & 172.0 & 177.0 & 183.1 & 291.0 & 299.3 & 311.9 \\
      \ac{RS3} & VaR  & 170.3 & 175.9 & 185.1 & 196.1 & 291.9 & 302.5 & 320.1 \\
      & CVaR                   & 170.7 & 177.6 & 189.3 & 203.7 & 292.3 & 304.2 & 326.9 \\
      \hline
      & Mean          & 171.1 & 173.4 & 178.4 & 185.2 & 302.1 & 430.4 & 591.3 \\
      \ac{SDPG} & VaR & 172.4 & 177.9 & 188.3 & 201.1 & 302.6 & 607.9 & 650.8 \\
      & CVaR          & 172.9 & 179.5 & 191.6 & 206.8 & 302.8 & 617.8 & 659.6 \\ 
      \bottomrule
    \end{tabular}
  \end{center}
\end{table}

\begin{figure}[t]
    \centering
    \includegraphics[width=0.49\textwidth]{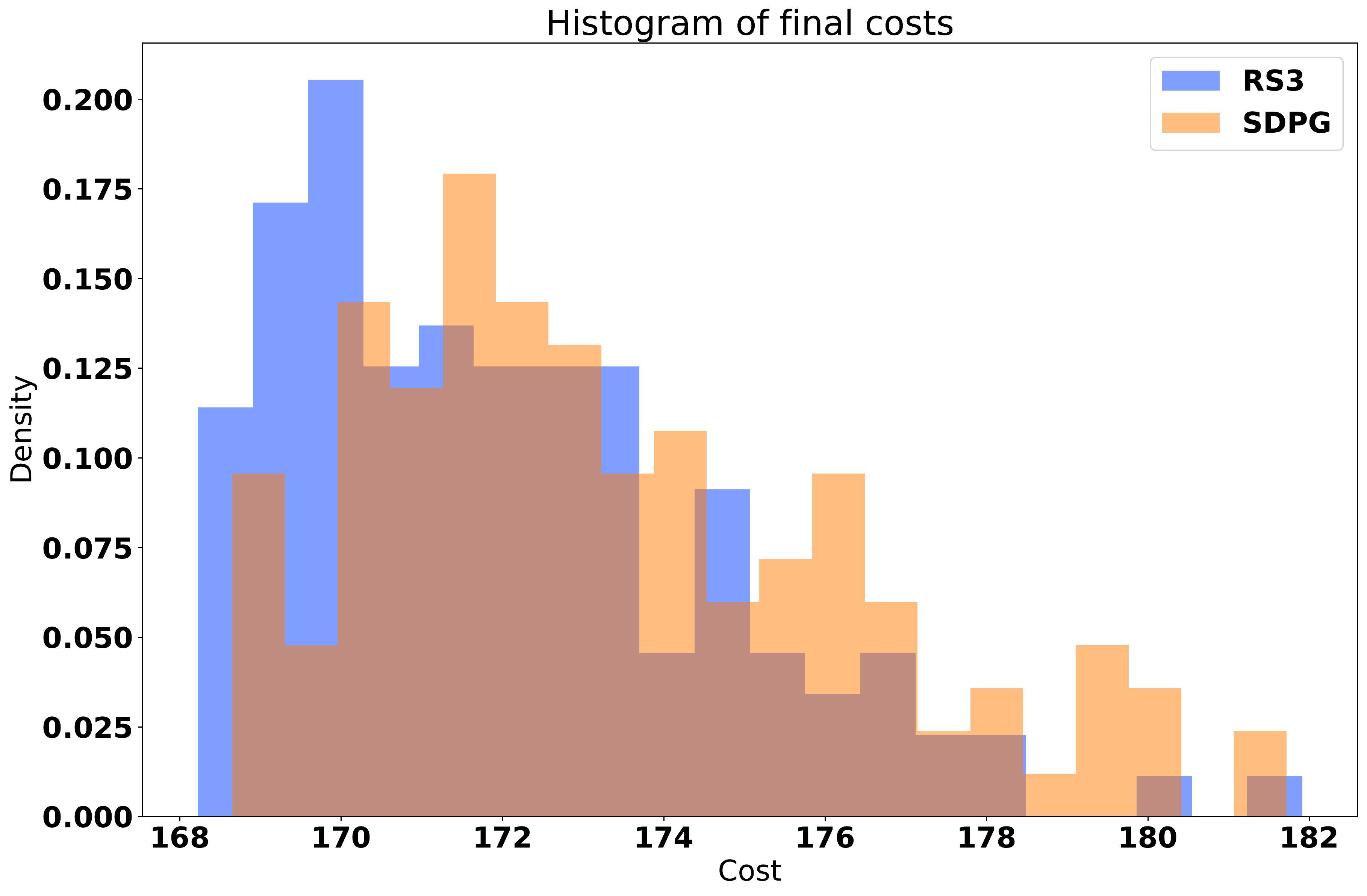}
    \includegraphics[width=0.49\textwidth]{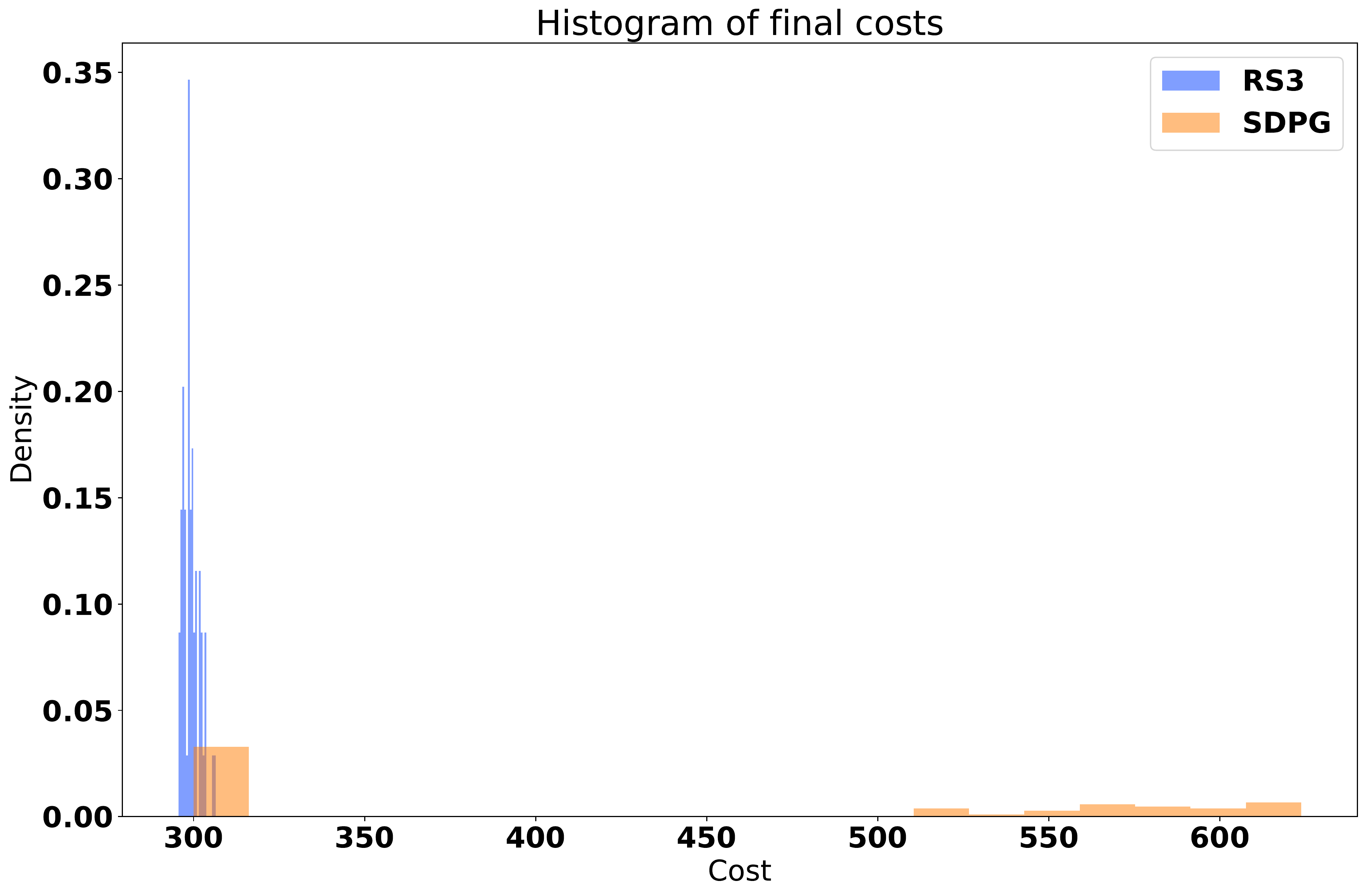}
    \vspace*{-5mm}
    \caption{The histograms of the final cost in the case of injected control noise sampled from $\mathcal{N}(0, 1).$ \textit{Left}: Pendulum, \textit{Right}: Cartpole. }
    \label{fig:comparison/SDPG/main_histogram}
\end{figure}

Typical \ac{RL} algorithms always receive some state feedback either fully or partially from environments. Thus, to fairly compare against \ac{SDPG}, we exploited an \ac{MPC} scheme in \ac{RS3} to implicitly receive the state feedback and perform a receding-horizon optimization.
In addition, as \ac{SDPG} is unable to handle uncertainty in the initial states and controls, we consider deterministic initial states and system dynamics with additive noise h the control channels.
To match the \ac{RS3} framework, the Gym environment's controls were modified to be continuous and use a quadratic cost function instead of the typical \ac{RL} reward function -1, 0, or +1 implemented in Gym.
The cost function used in the simulation can be found in 
Appendix \ref{sec:cost_function}.
All other training parameters for risk-sensitive \ac{SDPG} were the same as the parameters used in the original work \cite{singh2020cvarsdpg}.

\begin{wrapfigure}{r}{0.5\textwidth}
    \centering
    \includegraphics[width=0.5\textwidth]{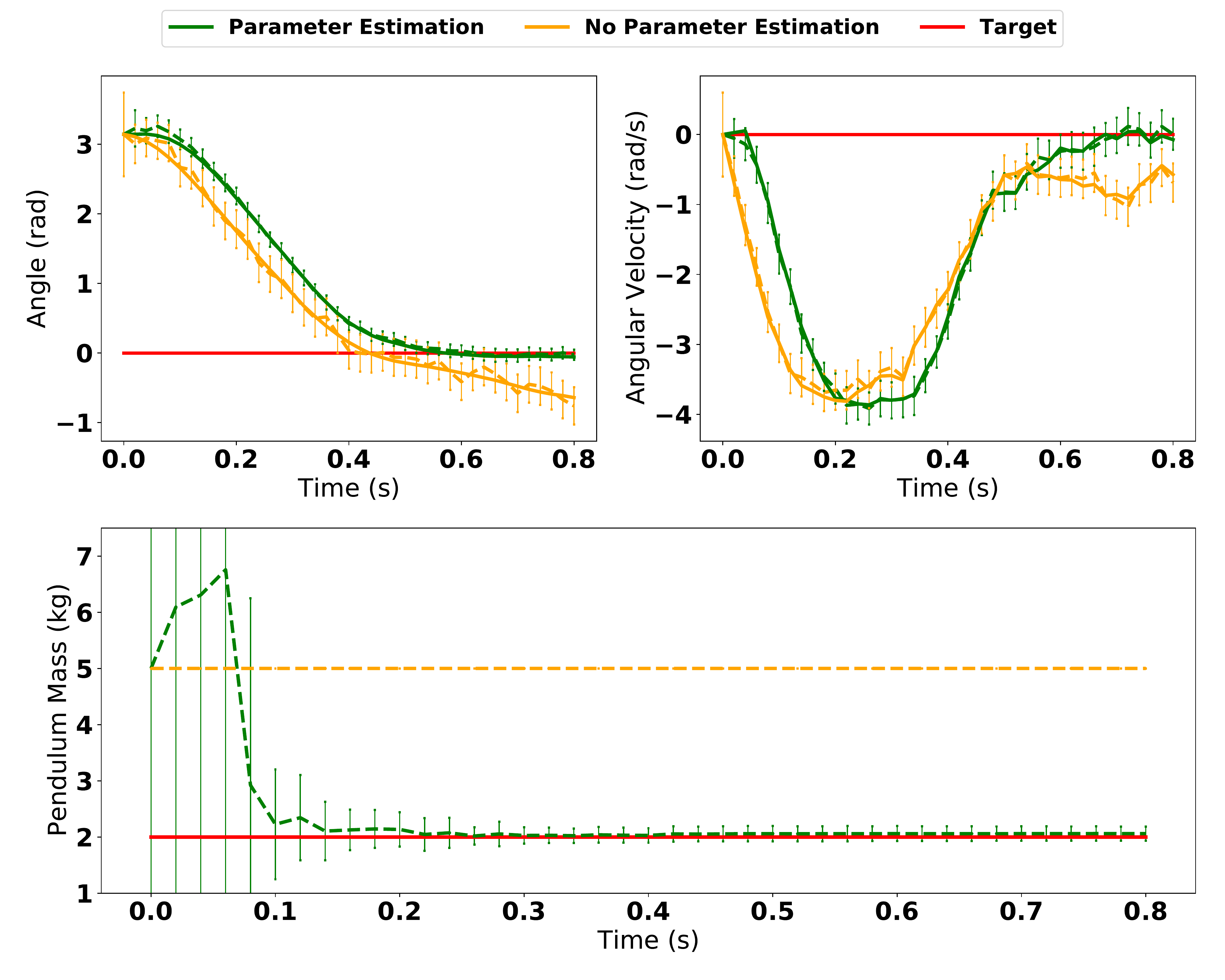}
    \vspace*{-10 mm}
    \caption{Nonlinear belief space optimization with uncertain pendulum mass in the Pendulum problem.}
    \vspace{-15 pt}  
    \label{fig:parameter_estimation/pendulum}
\end{wrapfigure}

Under the aforementioned conditions, \ac{RS3} is shown to outperform \ac{SDPG} overall by converging to a lower \ac{CVaR} value, especially in the case of larger noise levels. 
The mean, \ac{VaR}, and \ac{CVaR} values of the final costs obtained from both algorithms for the pendulum and cartpole simulation are shown in \Cref{table:comparison_sdpg} and a comparison of the histogram of the final costs are shown in \Cref{fig:comparison/SDPG/main_histogram}. The state, control, and cost histograms for all the simulation in \Cref{table:comparison_sdpg} can be found in 
\Cref{sec:appendix_comparison_sdpg}.

It is clearly shown in \Cref{fig:comparison/SDPG/main_histogram} that the distribution of the final cost has sharper tail on the high cost region in \ac{RS3}'s results compared to \ac{SDPG}'s. As a result, the mean, \ac{VaR}, and \ac{CVaR} of \ac{RS3}'s final costs are smaller than \ac{SDPG}'s.

The reason why our method outperforms the \ac{RL} framework is that we perform online update of our policy whereas the \ac{RL} policy is fixed after training.
This disadvantage of \ac{RL} algorithms comes from the nature of \ac{RL}. Once a model is trained on a specific dataset or with a specific noise profile, the model fails to output correct predictions under a new environment or given unseen inputs or noise. Our online optimization scheme solves this issue and fits better in risk-sensitive control.



\subsection{Belief Space Optimization}

\begin{wrapfigure}{h}{0.5\textwidth}
\vspace{-12 pt}
\centering
\includegraphics[width=0.5\textwidth]{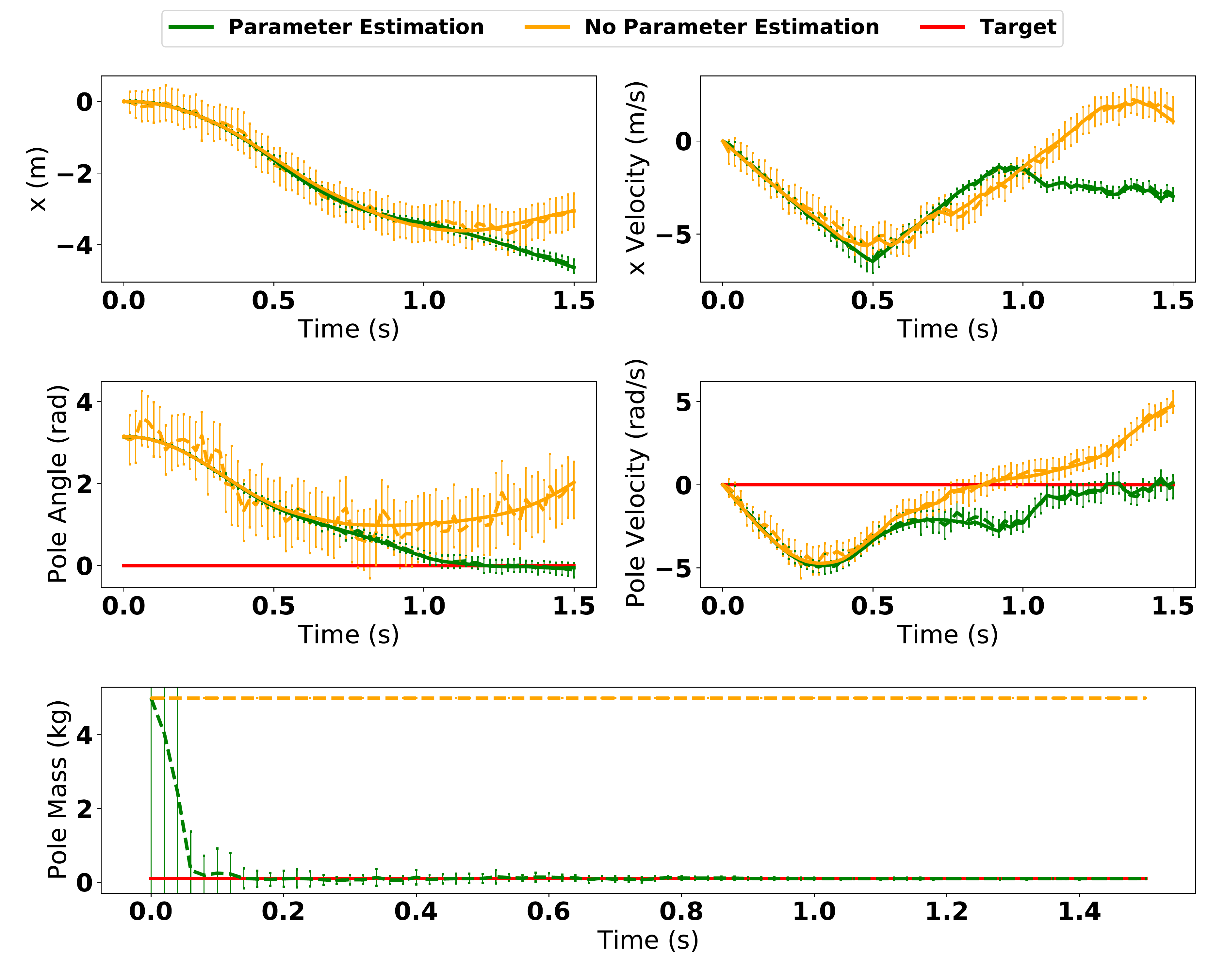}
    \vspace*{-10mm}
    \caption{Nonlinear belief space optimization with uncertain pole mass in the Cartpole problem.}
  \vspace{-15 pt}  
\label{fig:parameter_estimation/cartpole}
\end{wrapfigure}

We next show results for the uncertain parameter case from the pendulum, cartpole and quadcopter systems. In each trajectory plot, the dotted lines represent estimates from the particle filter with the error bars showing the $\pm3\sigma$ uncertainties of the nonlinear belief. The solid line represents the ground truth states.

\textbf{Pendulum: }
We first apply \ac{RS3} to a pendulum for a swingup task with unknown pendulum mass.
We assume deterministic initial condition and state transition model.
The pendulum's true mass is set to $2$ kg.
The prior for pendulum mass is set to be $\calN(5.0, 4.0)$.
The initial states $\mathbf{x} = [\theta, \dot{\theta}]$ are drawn from a normal distribution with mean $[\pi, 0]$ and covariance matrix $\textrm{diag}([0.1, 0.1])$. We assume full-state observability with additive measurement noise $\xi \sim \mathcal{N}(0, 1)$.
From \Cref{fig:parameter_estimation/pendulum}, we can observe that \ac{RS3} is able to correctly estimate the mass of the pendulum in the parameter estimation case. Without parameter estimation, \ac{RS3} overestimates the control effort required and overshoots the target angle.

\textbf{Cartpole: }
We apply the proposed algorithm to the task of cartpole swingup with with unknown pole mass.
The prior over the mass of the pole is a normal distribution $\mathcal{N}(5.0, 5.0)$ and the true value is $\SI{0.1}{kg}$.
Our algorithm is able to learn the true mass of the pole and successfully perform a swing up (\Cref{fig:parameter_estimation/cartpole}).
We compare this with the case of not estimating the mass of the pole, where
the algorithm does not sample from the correct dynamics and is unable to correctly optimize for a trajectory that successfully swings up.

\begin{figure}[t]
    \centering
    \includegraphics[width=0.9\textwidth]{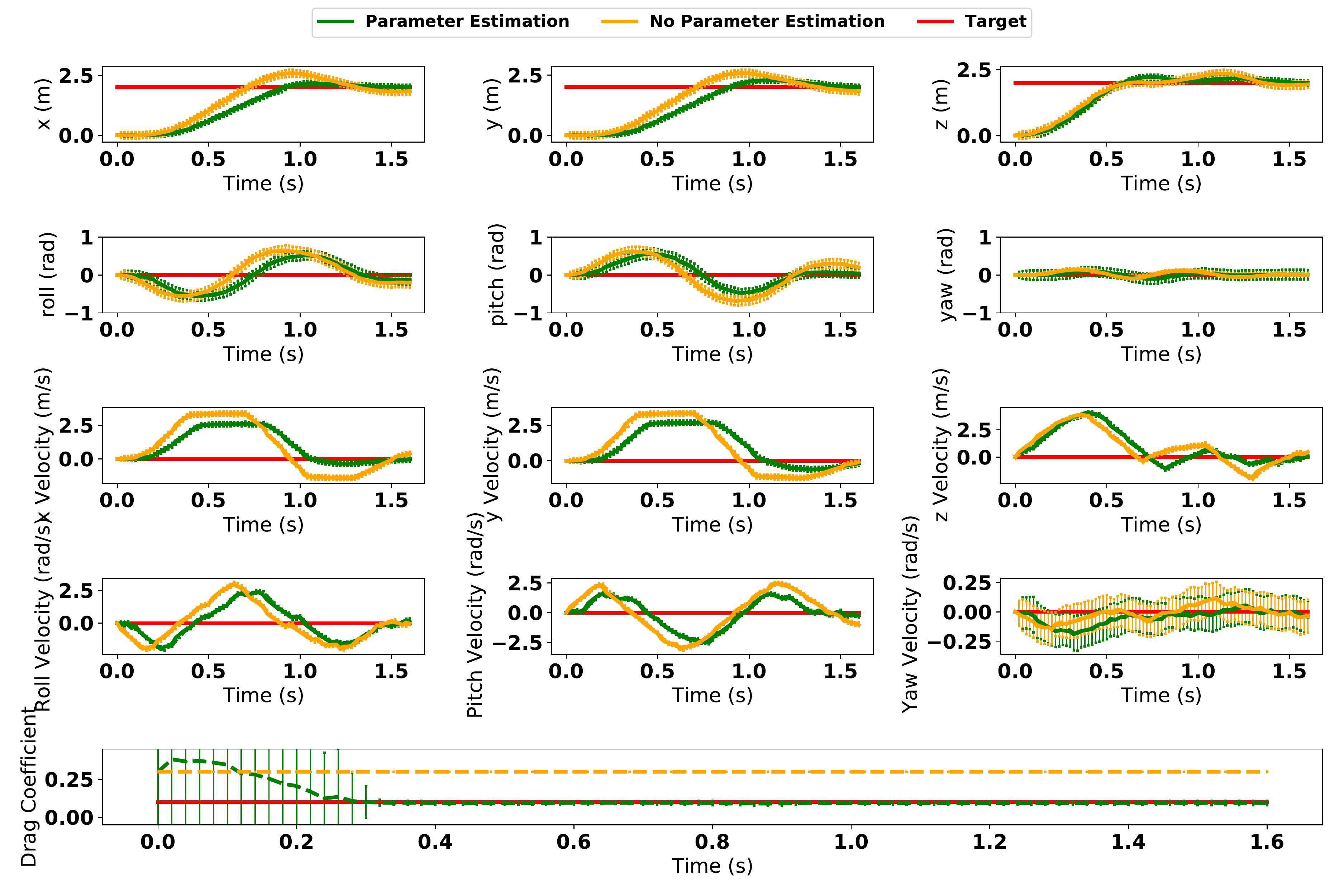}
    \vspace*{-5mm}
    \caption{Nonlinear belief space optimization with uncertain drag coefficient in the Quadcopter problem.
    \ac{RS3} with parameter estimation is able to converge much closer in roll, pitch and pitch velocity compared to without parameter estimation.}
    \label{fig:parameter_estimation/quadcopter}
\end{figure}

\textbf{Quadcopter: }
Finally, we apply our algorithm to the quadcopter system (dynamics can be found in \cite{elkholy2014dynamic}), where the task is to fly a quadcopter with states $[x, y, z, \dot{x}, \dot{y}, \dot{z}, r, p, y, \dot{r}, \dot{p}, \dot{y}]$ from $[0, 0, 0]$ to $[2, 2, 2]$.
The drag coefficient of the system is unknown, the prior over the drag is a normal distribution $\mathcal{N}(0.5, 0.5)$ and the true value is $0.1$.
The algorithm is once again able to learn the correct drag coefficient and manages to pilot the quadcopter to the target position without significant overshoot despite the drag coefficient being $500\%$ larger than the mean of the prior (\Cref{fig:parameter_estimation/quadcopter}). For the case where we do not perform parameter estimation, since the prior of the drag coefficient is greater than the actual value, the control policy found by the \ac{RS3} framework results in overshooting behavior before convergence, although it still manages to converge to the target state due to the robustness from optimizing for \ac{CVaR}. 

\section{Conclusion}
\label{sec:conclusion}

In this paper we introduced a general framework for \ac{CVaR} optimization for dynamical systems. The resulting algorithm, \ac{RS3}, is capable of handling uncertainties arising from uncertain initial conditions, unknown model parameters and system stochasticity. The algorithm can be readily combined with any filter for belief space risk sensitive control. We compared \ac{RS3} against \ac{SDPG} on the systems of a pendulum and cartpole and demonstrated outperformance in terms of final \ac{CVaR} cost. In addition, we combined \ac{RS3} with a particle filter for adaptive risk-sensitive control on non-Gaussian belief under different sources of uncertainty.

\section*{Acknowledgement}
This work is supported by NASA Langley and the NSF-CPS award \#1932288.

\newpage
\bibliography{references}
\bibliographystyle{unsrt}

\clearpage

\appendix

\section{Coherent Risk Measures}
\label{sec:coherency}
A risk measure $\rho:\mathcal{L}\rightarrow \mathbb{R}\cup\{+\infty\}$ is coherent if it satisfies certain simple mathematical properties including:
\begin{enumerate}
    \item \textbf{Monotonicity:} If $X_1, X_2\in\mathcal{L}$ and $X_1 \leq X_2$ a.s., then  $\rho(X_1) \leq \rho(X_2)$.
    \item \textbf{Sub-additivity:} If $X_1, X_2\in\mathcal{L}$, then $\rho(X_1 + X_2) \leq \rho(X_1) + \rho(X_2)$.
    \item \textbf{Positive homogeneity:} If $\alpha \geq 0$ and $X\in\mathcal{L}$, then $\rho(\alpha X)=\alpha\rho(X)$.
    \item \textbf{Translation invariance: }If $X\in\mathcal{L}$ and $c\in\mathbb{R}$ is constant, then $\rho(X+c)=\rho(X)+c$.
\end{enumerate}
Note that \ac{VaR} does not respect the sub-additivity property.

\section{Stochastic Search Update Law Derivation}
\label{sec:derivation_appendix}

The stochastic search update law \eqref{eq:gass_update_law} can be derived as
\begin{align*}
    \nabla_{\theta_t} l(\theta)&= \frac{\nabla_{\theta}\Eb[S(-\CVaR^\gamma [J(X,U)])]}{\Eb[S(-\CVaR^\gamma [J(X,U)])]}\\
    &= \frac{\nabla_{\theta}\int_{\Omega_\eta} S(-\CVaR^\gamma [J(X,U)])p(\eta)\rd \eta}{\Eb[S(-\CVaR^\gamma [J(X,U)])]}\\
    &= \frac{\int_{\Omega_\eta} S(-\CVaR^\gamma [J(X,U)])\nabla_{\theta}p(\eta)\rd \eta}{\Eb[S(-\CVaR^\gamma [J(X,U)])]}\\
    &= \frac{\int_{\Omega_\eta} S(-\CVaR^\gamma [J(X,U)])p(\eta)\nabla_{\theta}\ln p(\eta)\rd \eta}{\Eb[S(-\CVaR^\gamma [J(X,U)])]}\\
    &= \frac{\int_{\Omega_\eta} S(-\CVaR^\gamma [J(X,U)])p(\eta)\nabla_{\theta}\big(\sum_{t=0}^{T-1} p(\eta_t;\theta_t)\big)\rd \eta}{\Eb[S(-\CVaR^\gamma [J(X,U)])]}\\
    &= \frac{\Eb[S(-\CVaR^\gamma[J(X,U)])\nabla_{\theta}\big(\sum_{t=0}^{T-1}\ln p(\eta_t;\theta_t)\big)]}{\Eb[S(-\CVaR^\gamma [J(X,U)])]}.
\end{align*}

\section{Particle Filter and RS3 for Uncertain Model Parameters}
\label{sec:pf_implementation}

In the case of unknown parameters $\phi$ in the model, we use the particle filter to learn both the states and uncertain model distribution. We first augment our state space with those parameters as $y=[x, \phi]\T$ and formulate the new dynamics as
\begin{align*}
    y_{t+1} &= \begin{bmatrix}f(x_t,u_t,\phi_t) \\ 1 \end{bmatrix} + \begin{bmatrix} 0 \\ w_t
    \end{bmatrix}, \quad \phi_0\sim p(\phi_0)\\
    z_t &= h(x_t) + v_t
\end{align*}
where $w_t$ is the process noise added to mitigate the sample impoverishment problem of the particle filter. Starting with deterministic initial states $x_0$ and a prior distribution on the parameters $p(\phi_0)$. The steps of this adaptive risk sensitive control framework at each timestep are

\begin{enumerate}
    \item Sample from initial distribution: $\{\hat{y}_t^i\}_{i=1,\dots,M} \sim p(y_t)$
    \item Compute policy: $\eta_t = \text{RS3}(\{\hat{y}_t^i\}_{i=1,\dots,M})$
    \item Apply control: $x_{t+1} \leftarrow F(x_t, \pi_{\eta_t}(x_t),0;\phi)$
    \item Obtain observation: $z_{t+1} \leftarrow h(x_{t+1}, \xi_{t+1})$
    \item Propagate particles: $\{\hat{x}_{t+1}^i\}_{i=1,\dots,M} \leftarrow \hat{F}(u_t, \{\hat{x}_t^i\}_{i=1,\dots,M},0;\{\hat{\phi}_t^i\}_{i=1,\dots,M})$
    \item Update likelihood: $\{q_i\}_{i=1,\dots,M} = L(z_{t+1}, \{\hat{x}_{t+1}^i\}_{i=1,\dots,M})$
\end{enumerate}
where $F$ is the true system dynamics with real parameters $\phi$ and $\hat{F}$ is the model for particle filter and \ac{RS3}, and both do not contain stochasticities in the dynamics. $\xi$ is the observation noise.

\section{Additional Simulation Results}
\subsection{Comparison against SDPG}
\label{sec:appendix_comparison_sdpg}
\subsubsection{Cost function}
\label{sec:cost_function}
We define the cost function $J(X, U)$ as
\begin{align}
    J(X, U) &= g(x_{T}) + \sum_{t=t_0}^{T-1}l(x_t, u_t),
\end{align}
where $x_t \in \mathbb{R}^{n_x}$ is the state at time $t$, $u_t \in \mathbb{R}^{n_u}$ is the control at time $t$, $g$ is the state cost at the final time $T$, and $l$ is the running cost at time $t$.
In the simulations, we set $g=l$ and the running cost is composed of the quadratic state cost $(x_t-x_{target})\T Q (x_t-x_{target})$ and the quadratic control cost $u_t\T R u_t$.

The system dynamics are from OpenAI Gym \cite{brockman2016openai} and the pendulum state is defined as $[\theta, \dot{\theta}]$, the pendulum angle and the angular velocity.
The cartpole state is defined as $[x, \dot{x}, \theta, \dot{\theta}]$, the cart position, position velocity, pole angle, and the pole angular velocity.
The state and control cost used in pendulum dynamics are $Q=\text{diag}([3, 0.01])$, where $\text{diag}()$ represents a diagonal matrix, and $R=0.01$.
In cartpole, the costs are $Q=\text{diag}([0.01, 0.1, 1, 0.1])$ and $R=0.001$.

For Pendulum, the initial state is $[-\pi, 0]$ and target state is $[0, 0]$ and for the Cartpole, the initial state is $[0,0,-\pi,0]$ and target state is $[0, 0, 0, 0]$.

\subsubsection{Pendulum}
\begin{figure}[H]
    \centering
    \includegraphics[width=0.49\textwidth]{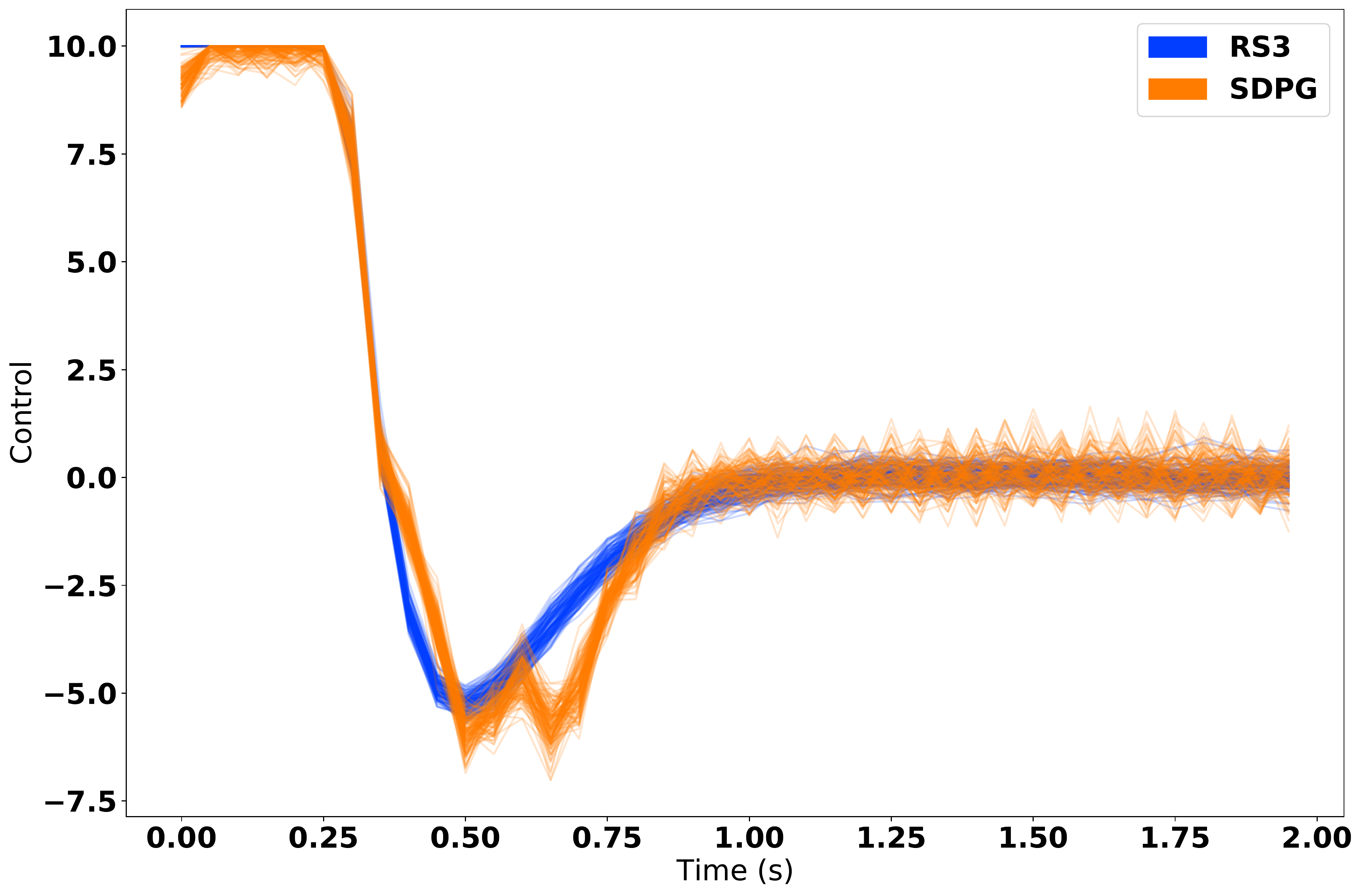}
    \includegraphics[width=0.49\textwidth]{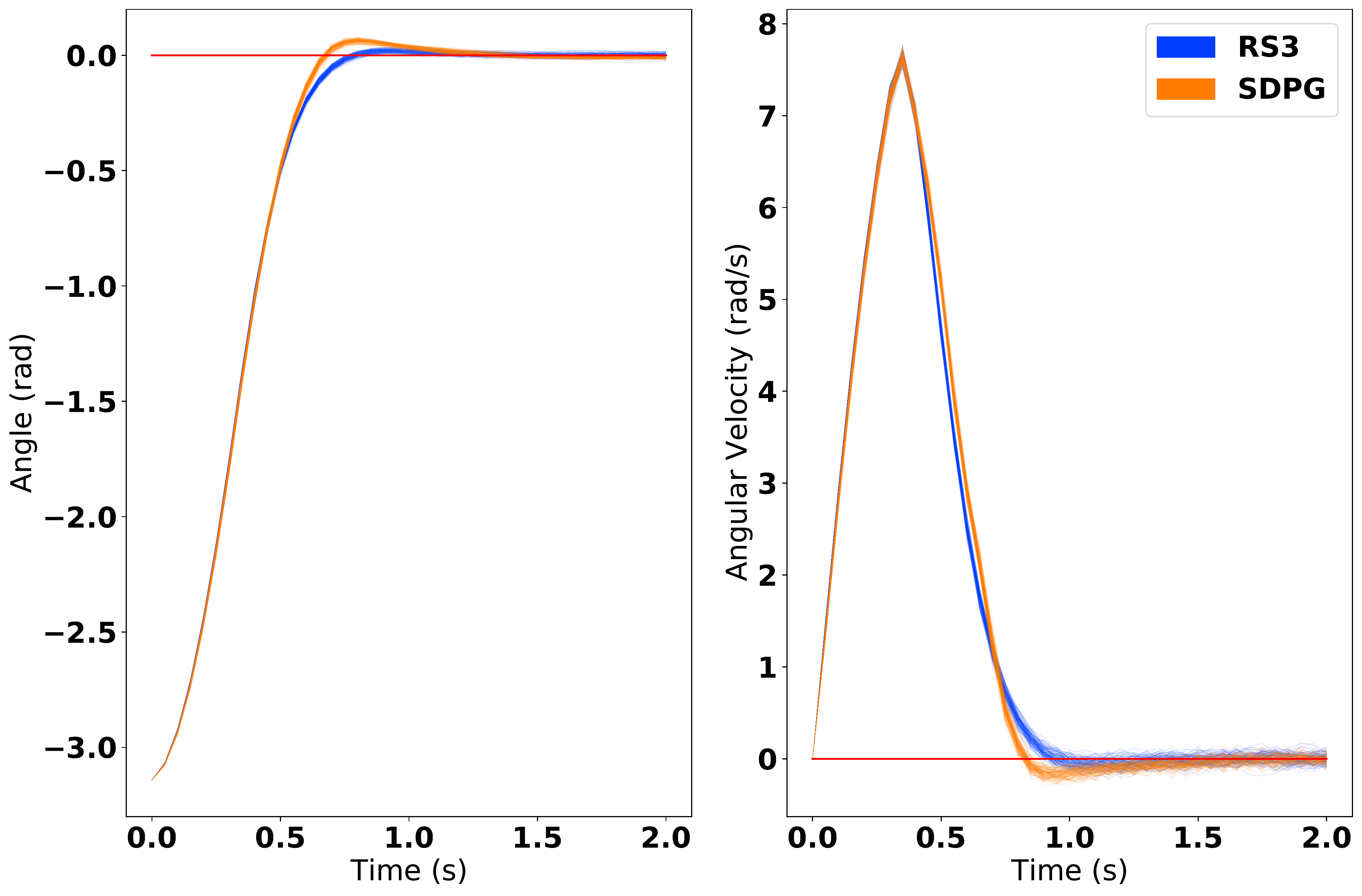}
    \includegraphics[width=0.49\textwidth]{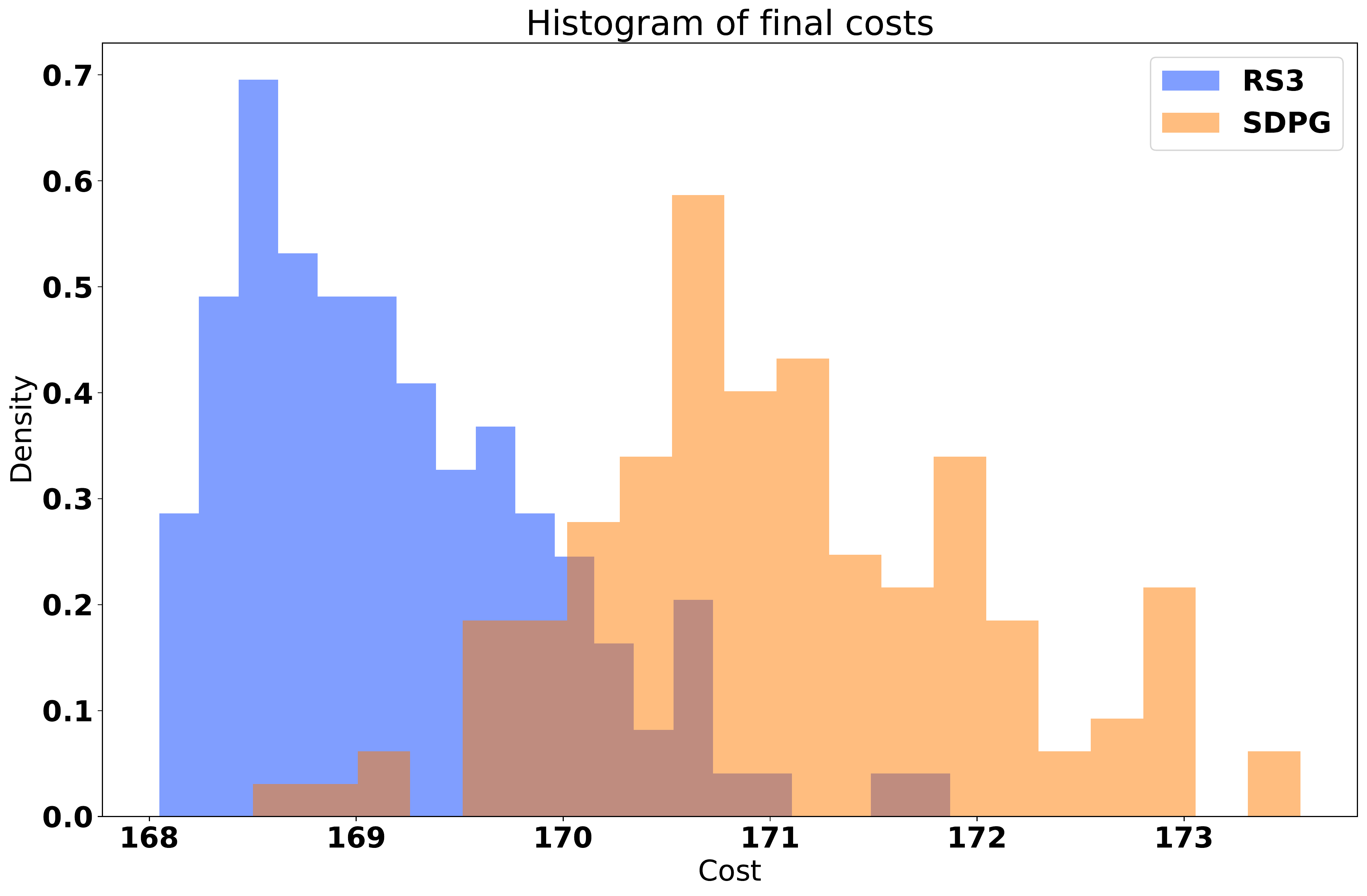}
    \includegraphics[width=0.49\textwidth]{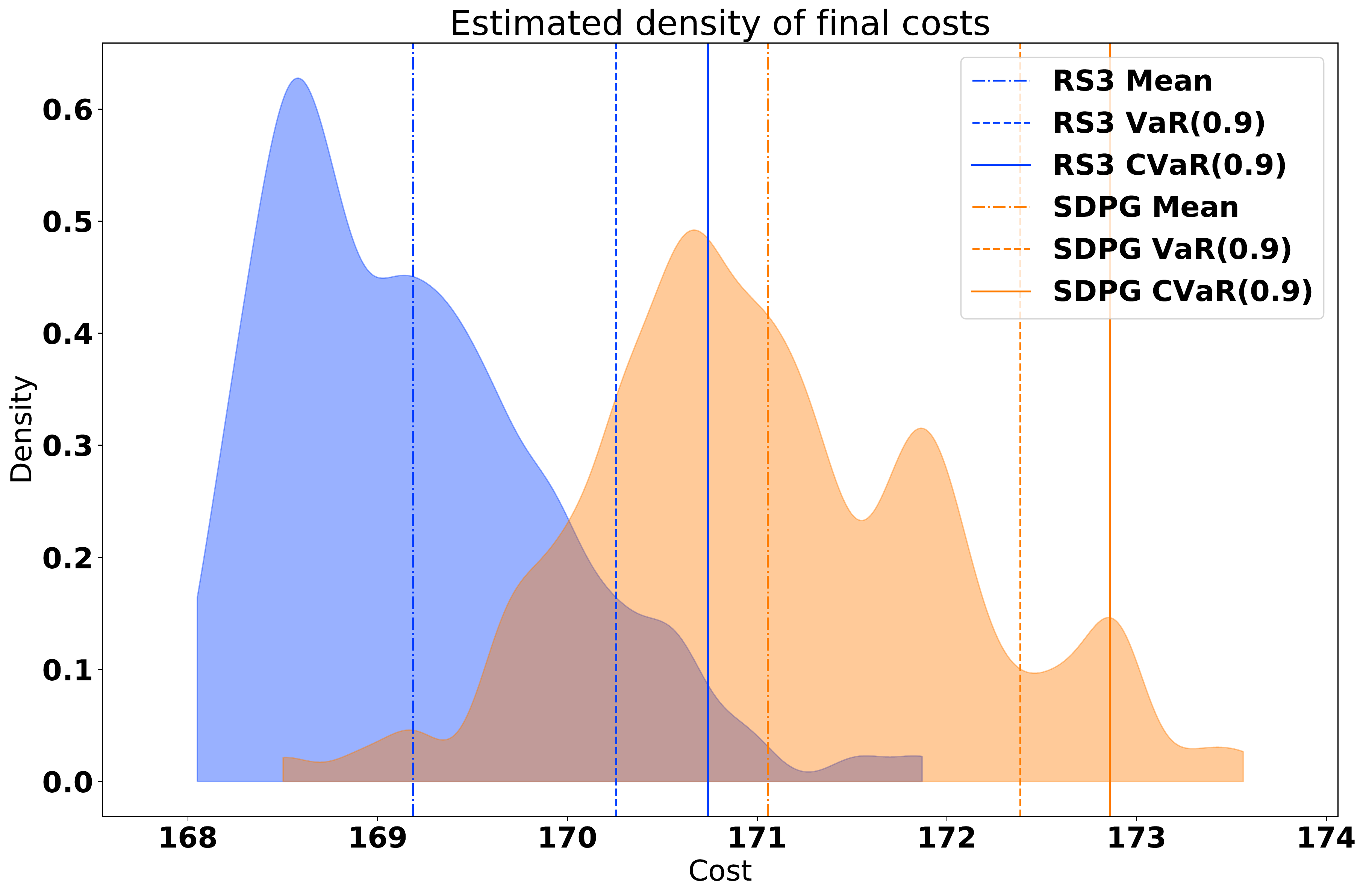}
    \caption{\textit{Top}: Control and state trajectories from the Pendulum problem with noise in the control channel. The control noise is $\mathcal{N}(0, 0.3^2).$ \textit{Bottom}: Cost histogram and the estimated p.d.f. of it. }
    \label{fig:comparison/SDPG/Pendulum/0.3}
\end{figure}

\begin{figure}[H]
    \centering
    \includegraphics[width=0.49\textwidth]{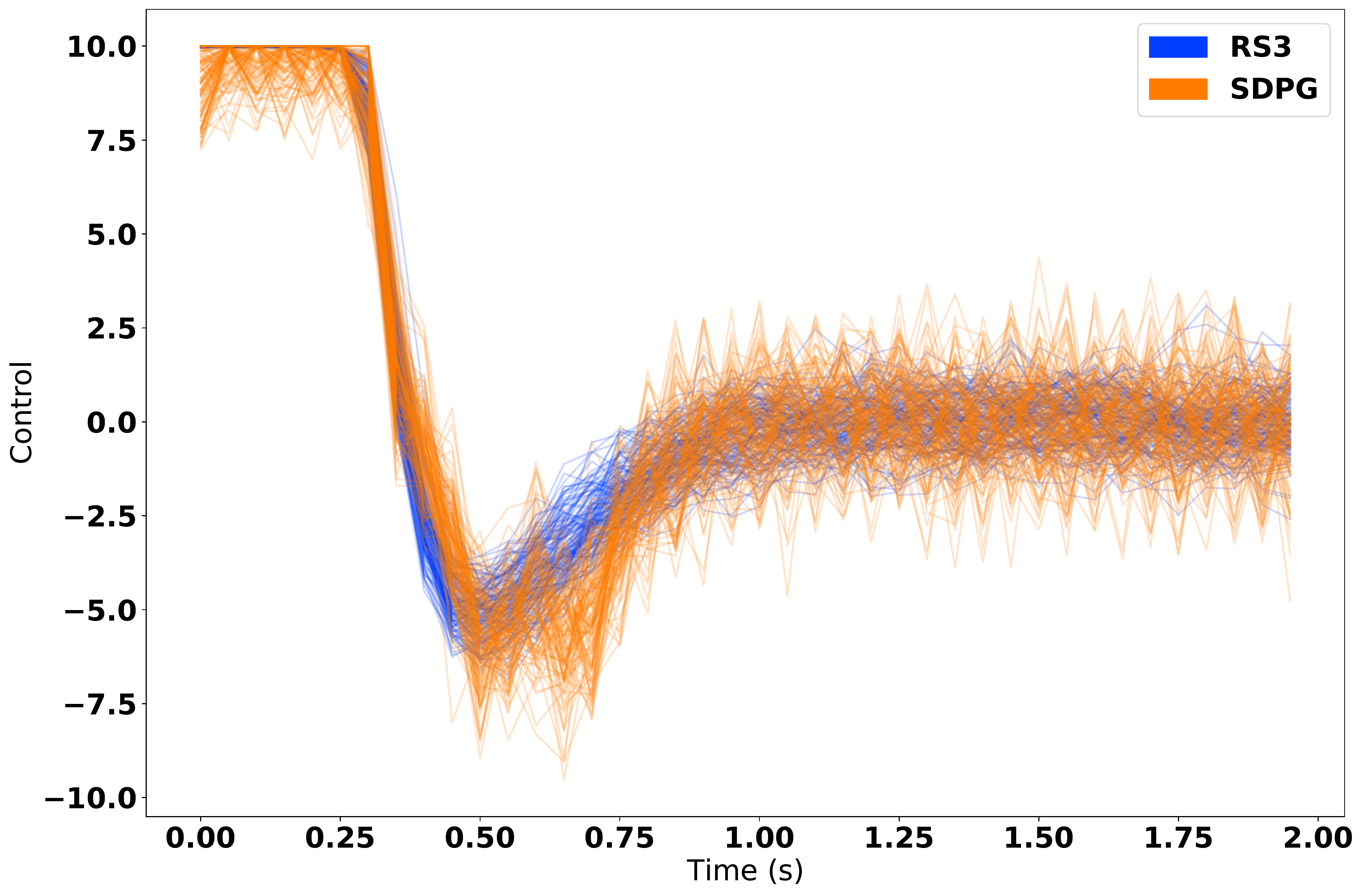}
    \includegraphics[width=0.49\textwidth]{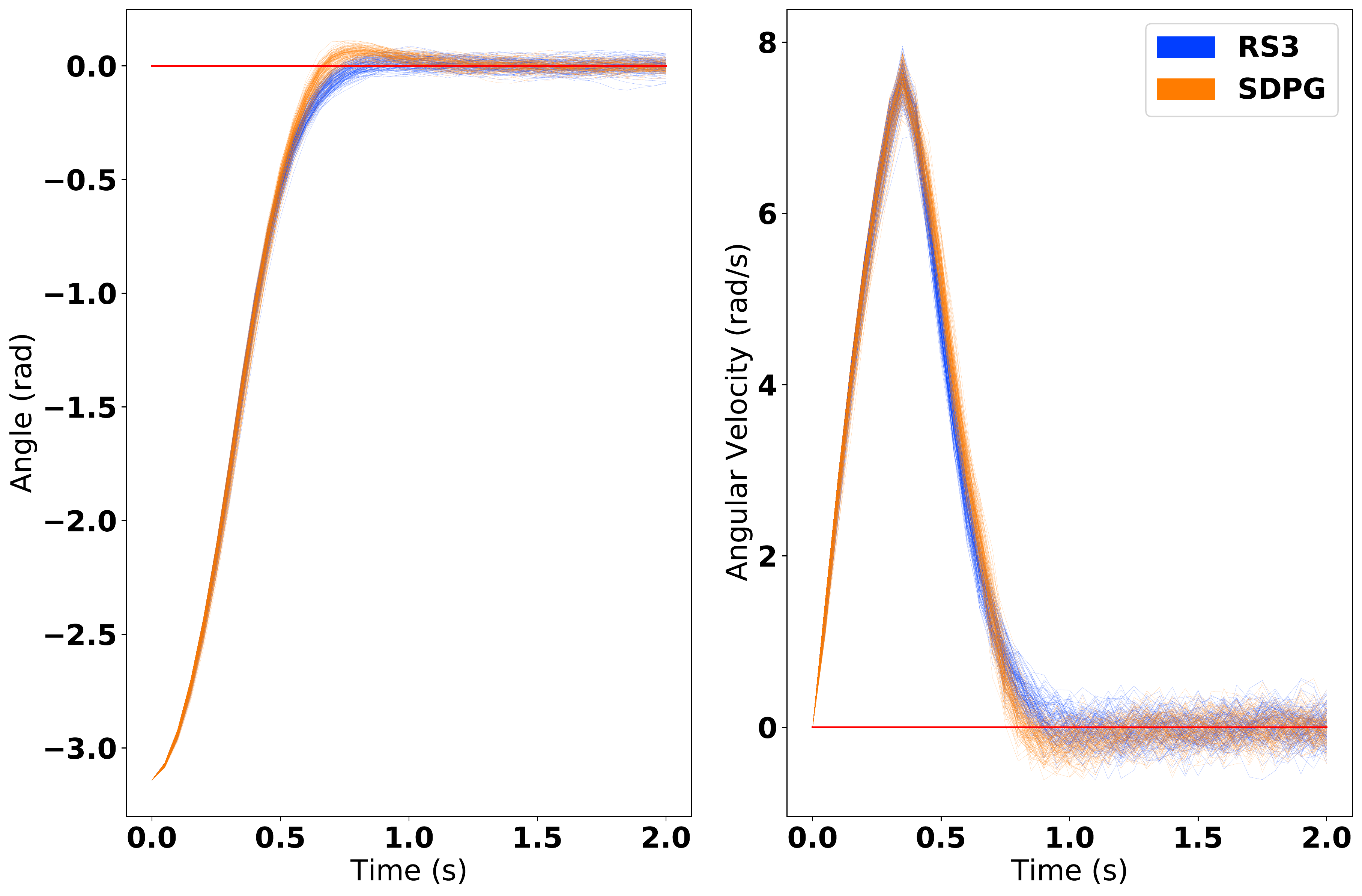}
    \includegraphics[width=0.49\textwidth]{figures/comparison/SDPG/Pendulum/1.0/cost_histogram.pdf}
    \includegraphics[width=0.49\textwidth]{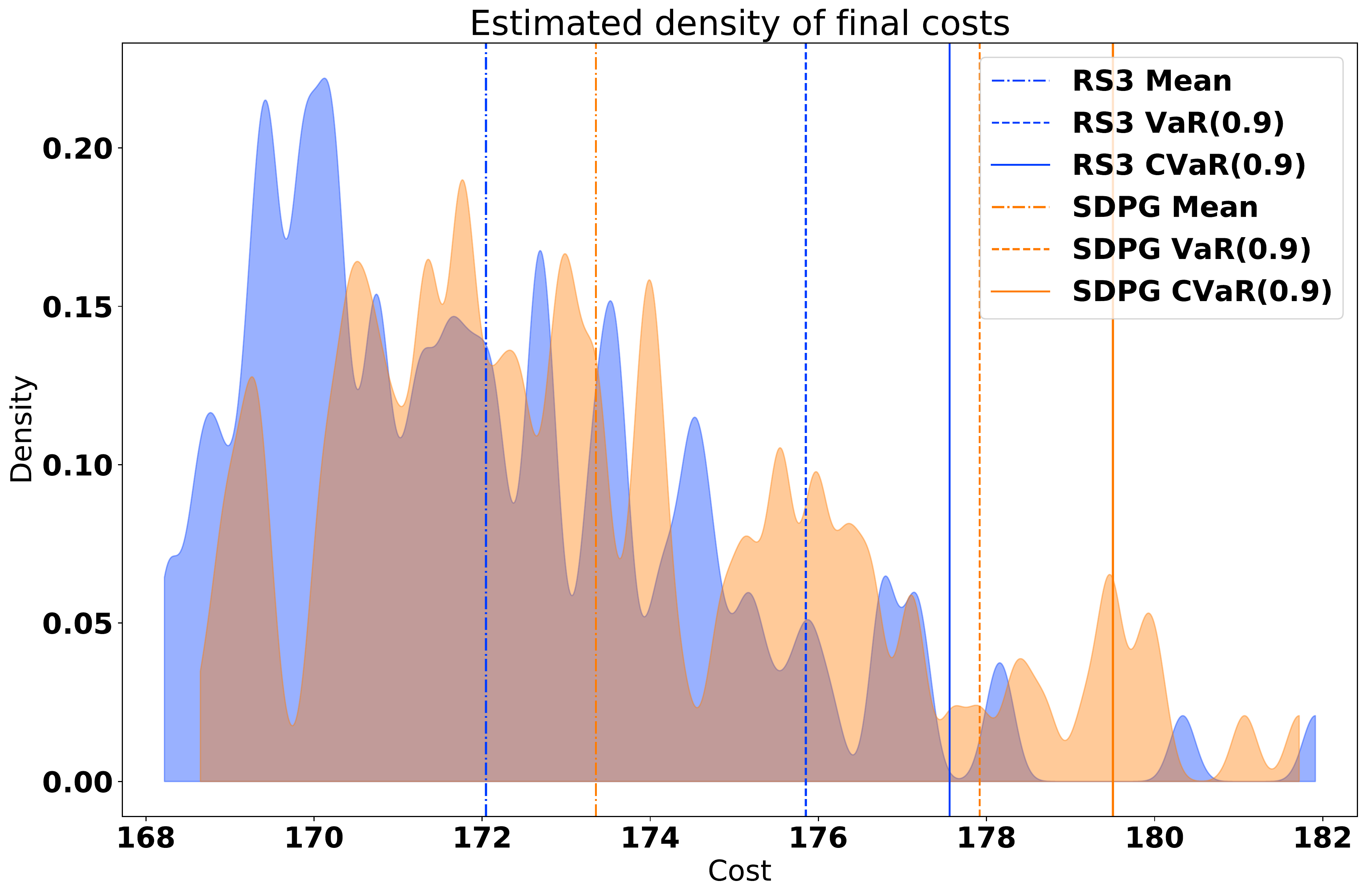}
    \caption{\textit{Top}: Control and state trajectories from the Pendulum problem with noise in the control channel. The control noise is $\mathcal{N}(0, 1.0^2).$ \textit{Bottom}: Cost histogram and the estimated p.d.f. of it. }
    \label{fig:comparison/SDPG/Pendulum/1.0}
\end{figure}

\begin{figure}[H]
    \centering
    \includegraphics[width=0.49\textwidth]{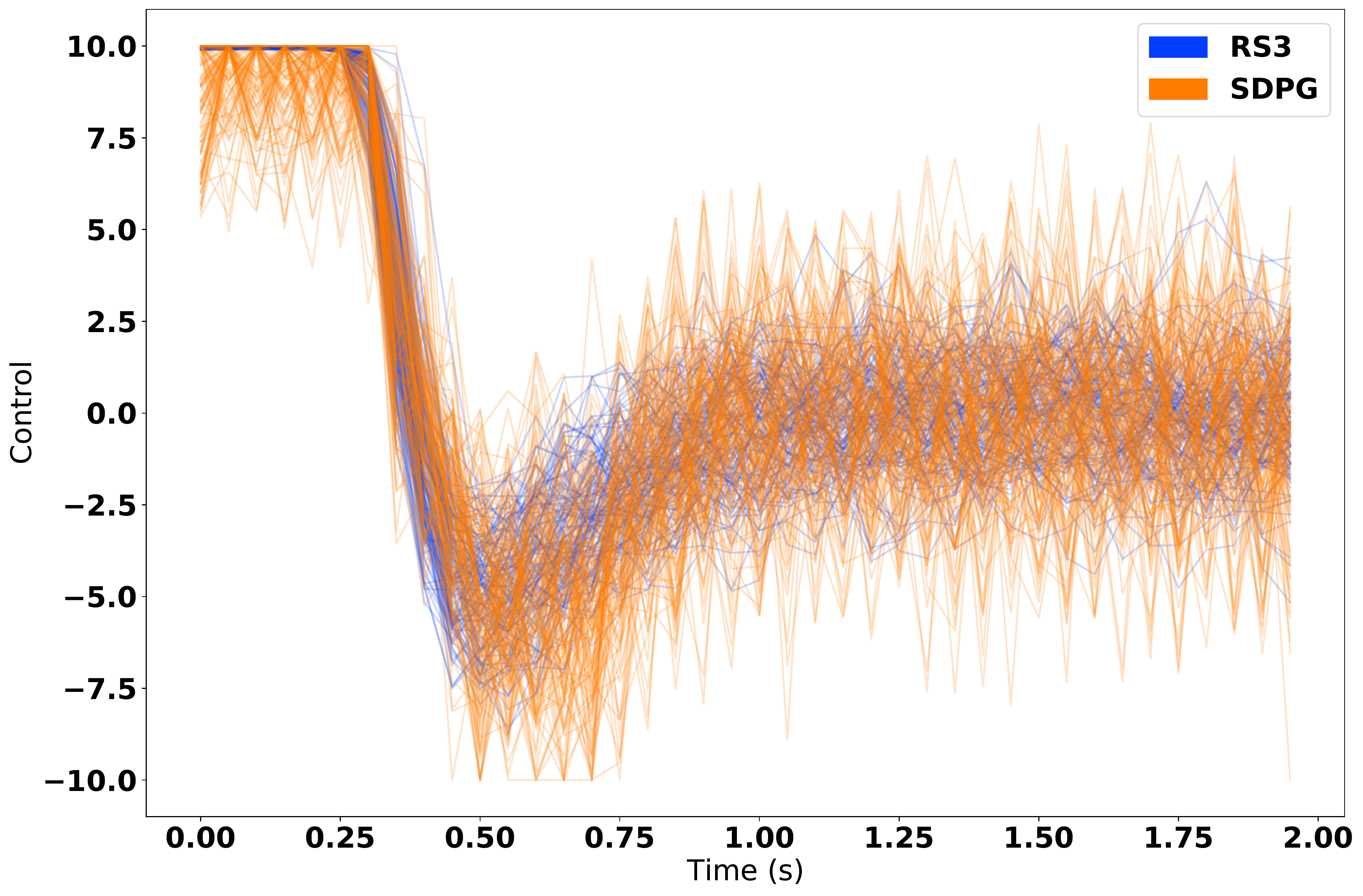}
    \includegraphics[width=0.49\textwidth]{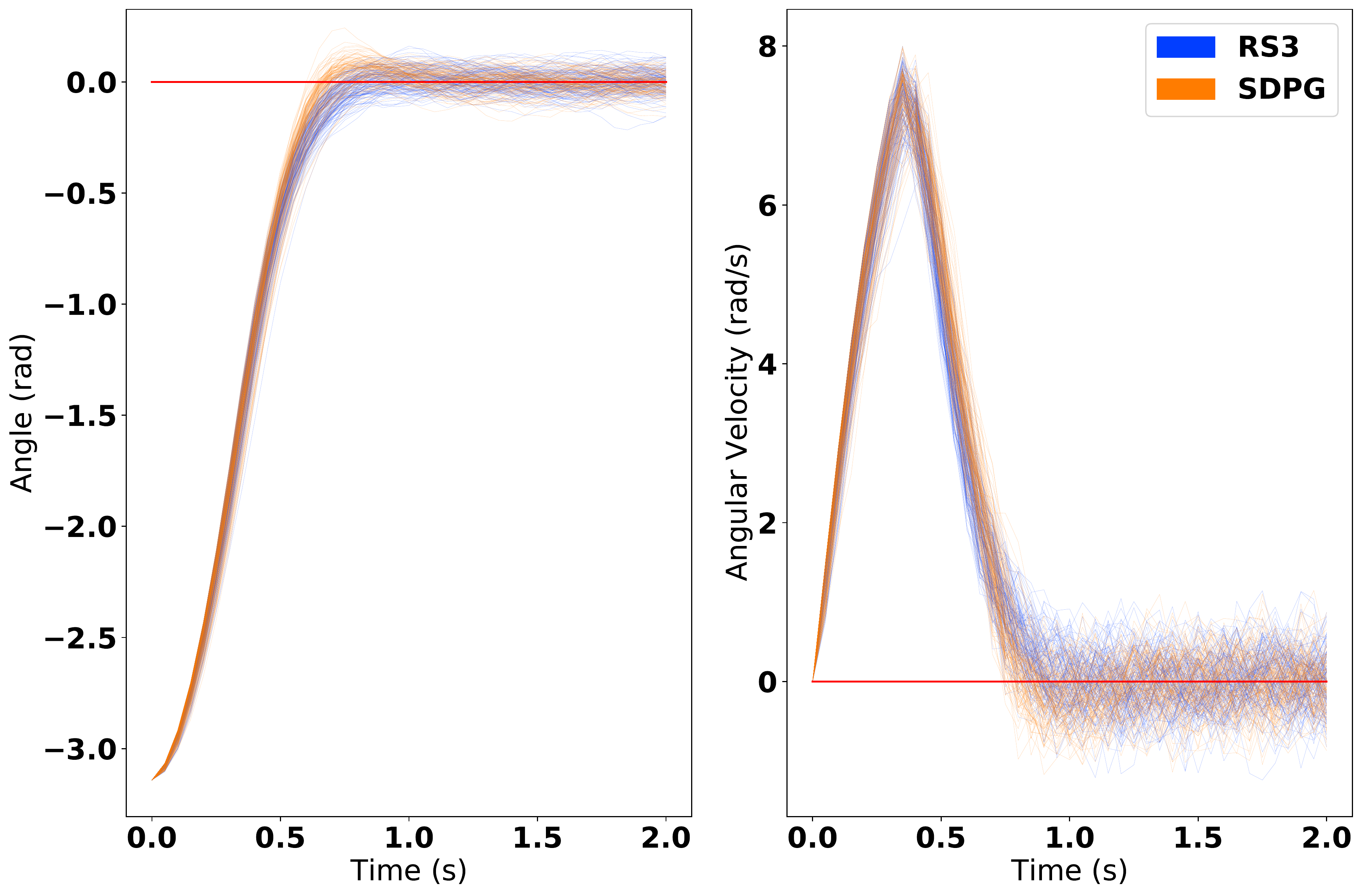}
    \includegraphics[width=0.49\textwidth]{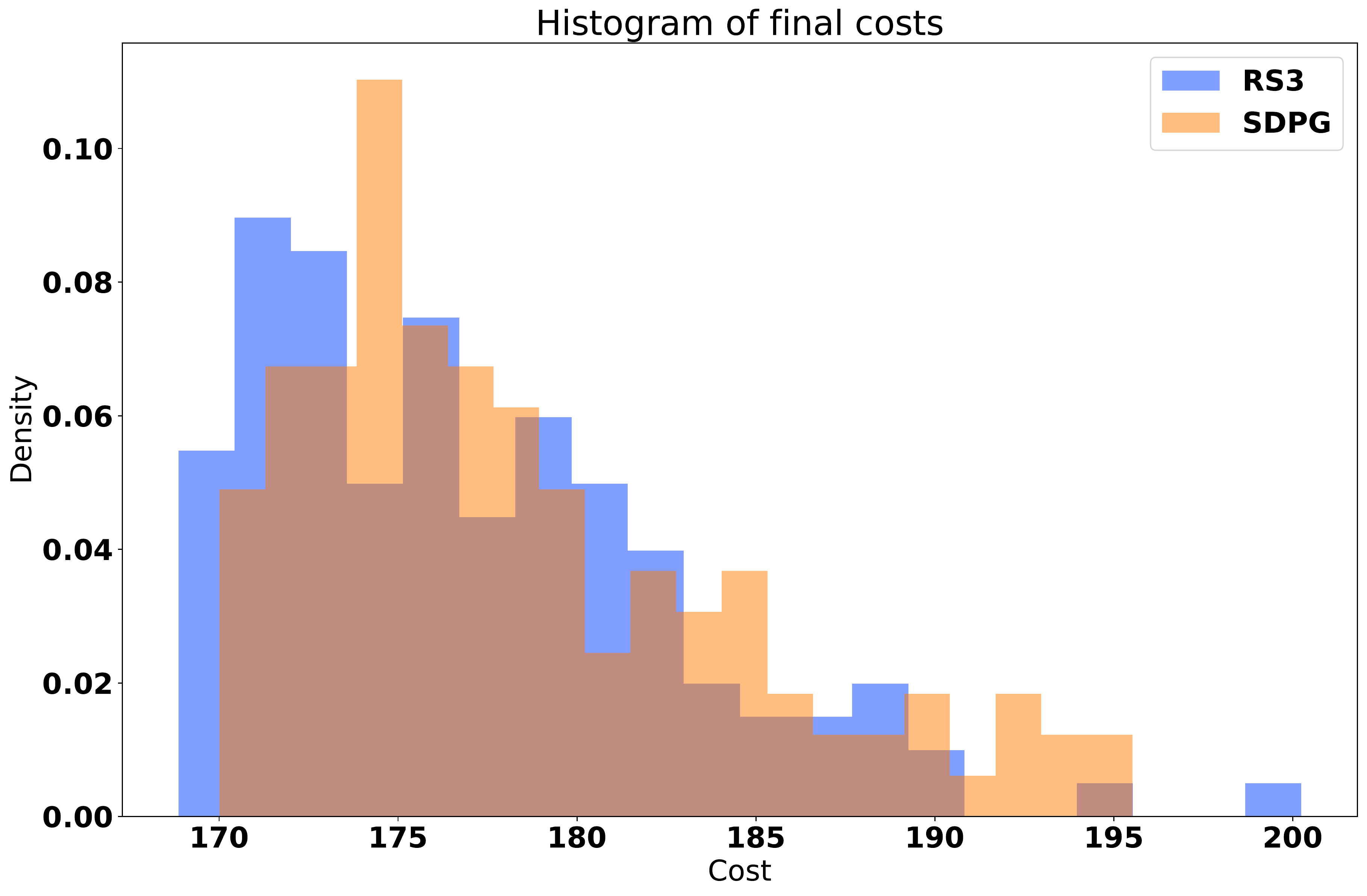}
    \includegraphics[width=0.49\textwidth]{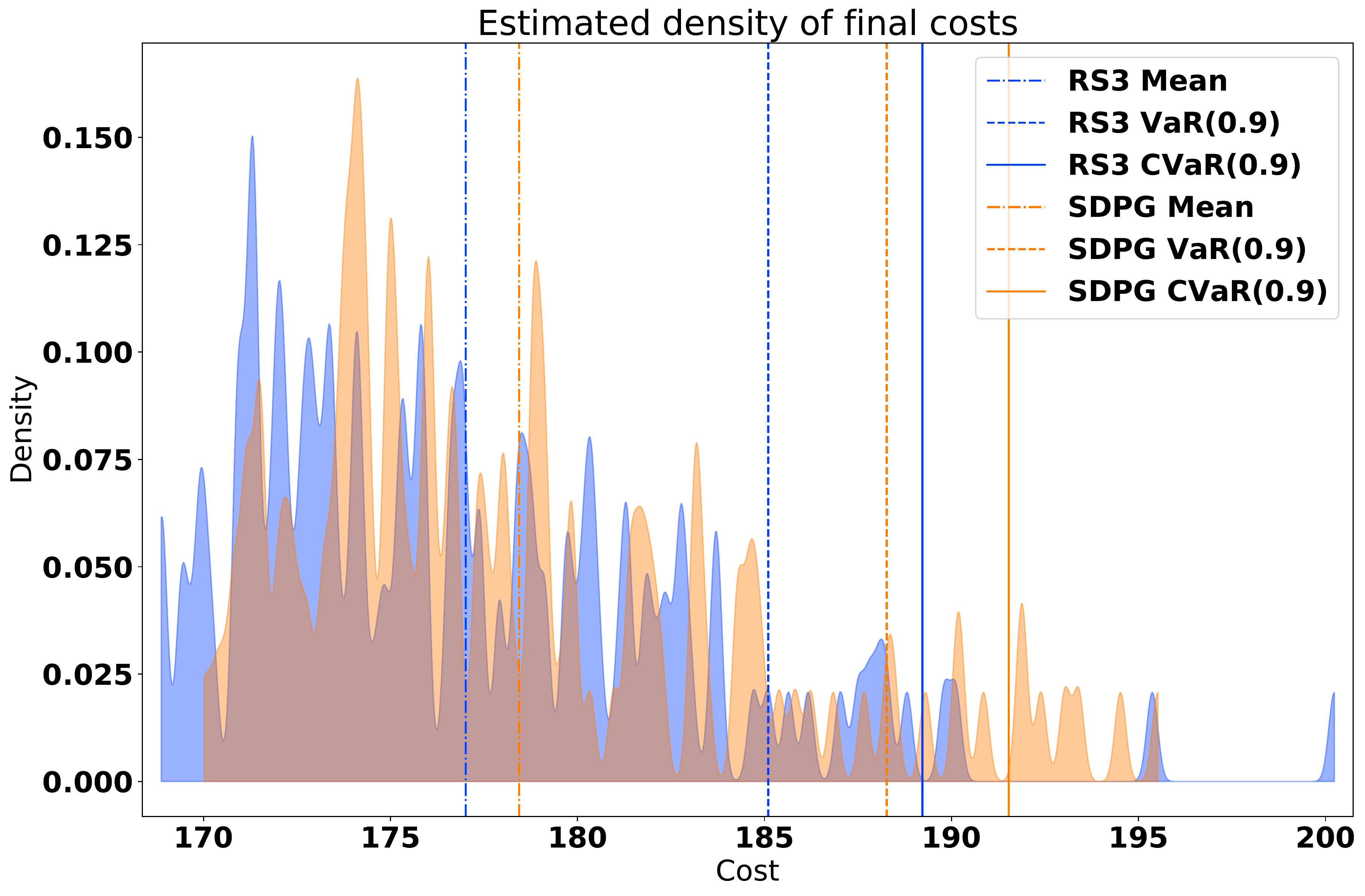}
    \caption{\textit{Top}: Control and state trajectories from the Pendulum problem with noise in the control channel. The control noise is $\mathcal{N}(0, 2.0^2).$ \textit{Bottom}: Cost histogram and the estimated p.d.f. of it. }
    \label{fig:comparison/SDPG/Pendulum/2.0}
\end{figure}

\begin{figure}[H]
    \centering
    \includegraphics[width=0.49\textwidth]{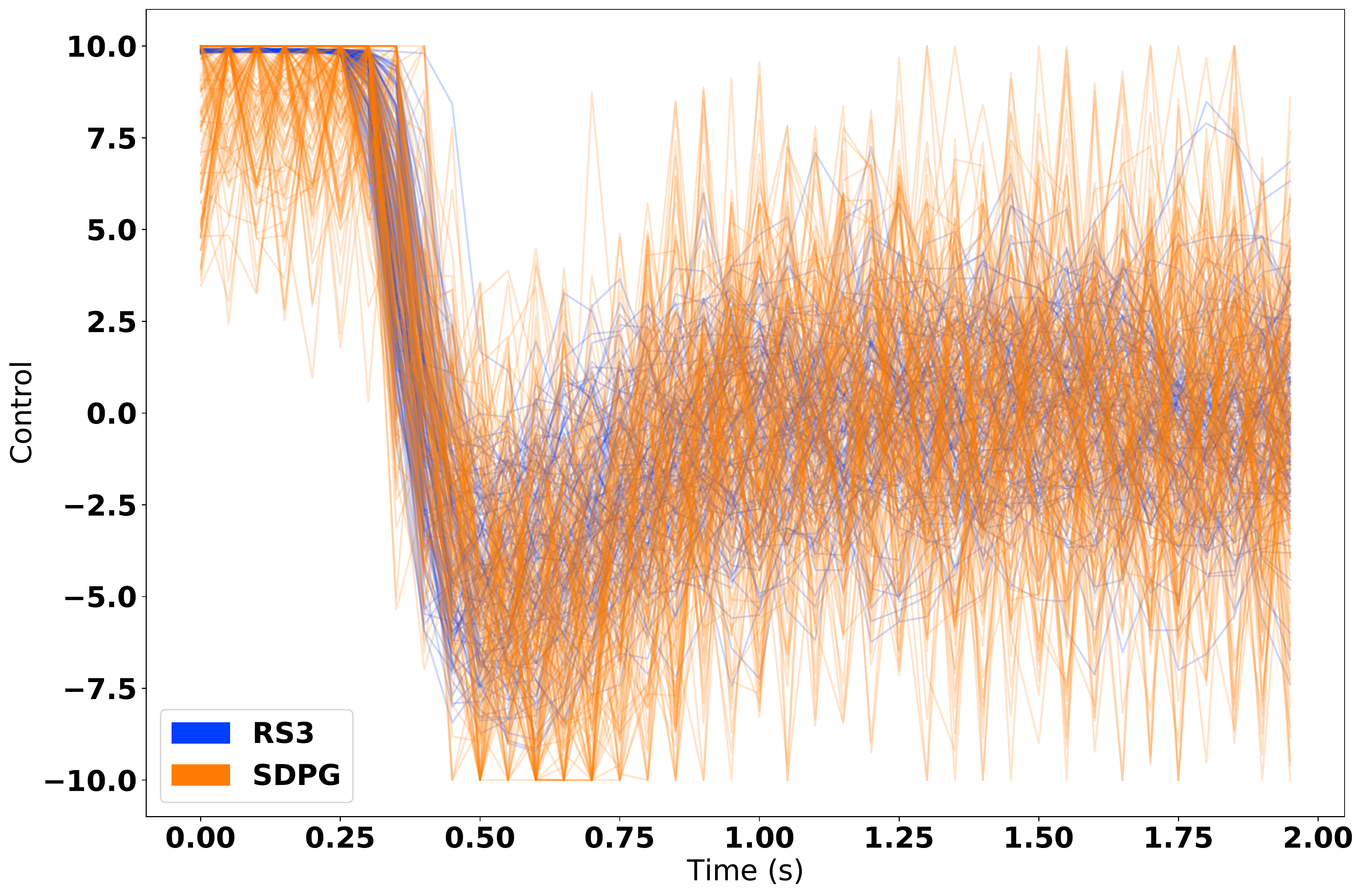}
    \includegraphics[width=0.49\textwidth]{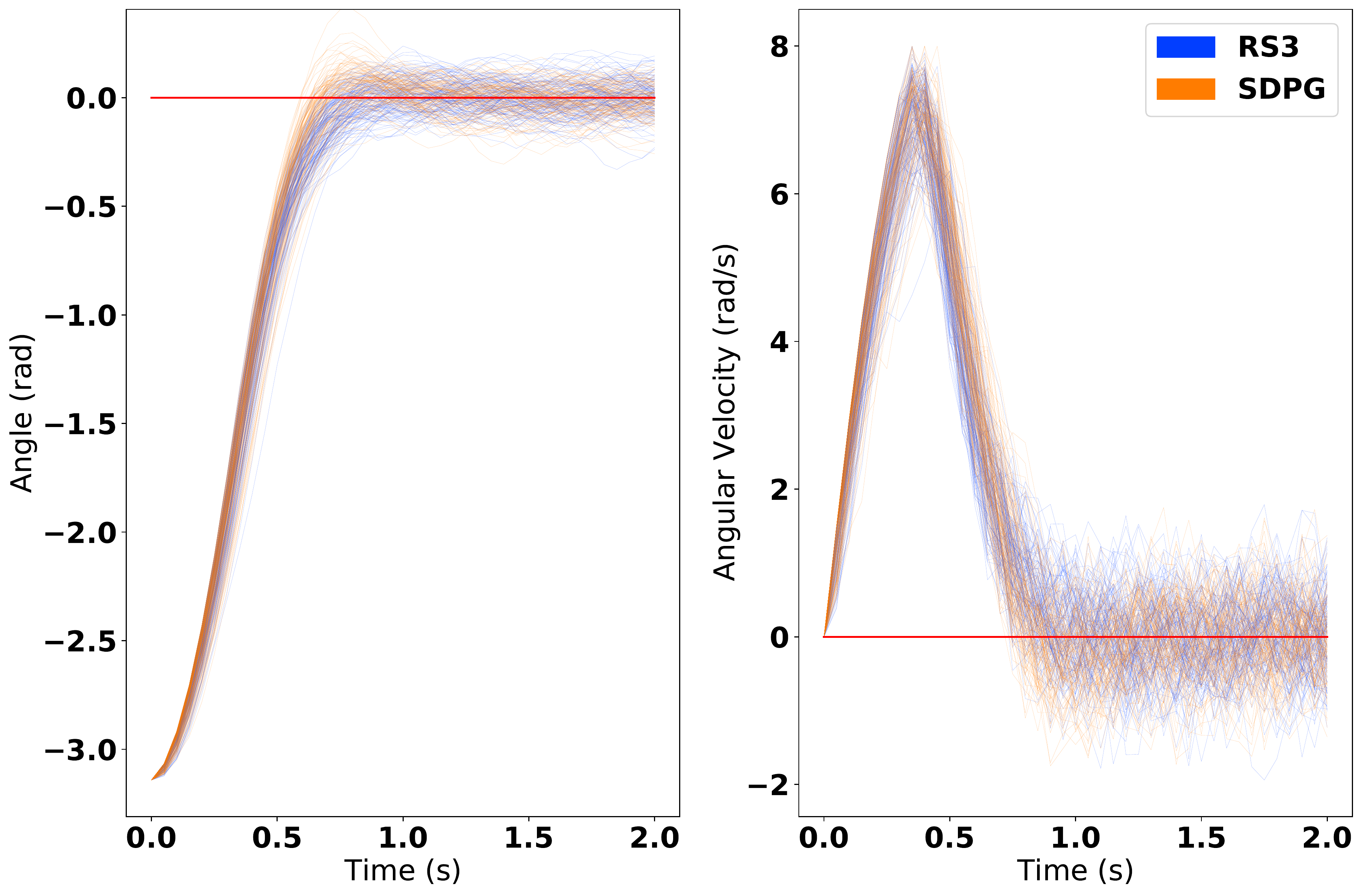}
    \includegraphics[width=0.49\textwidth]{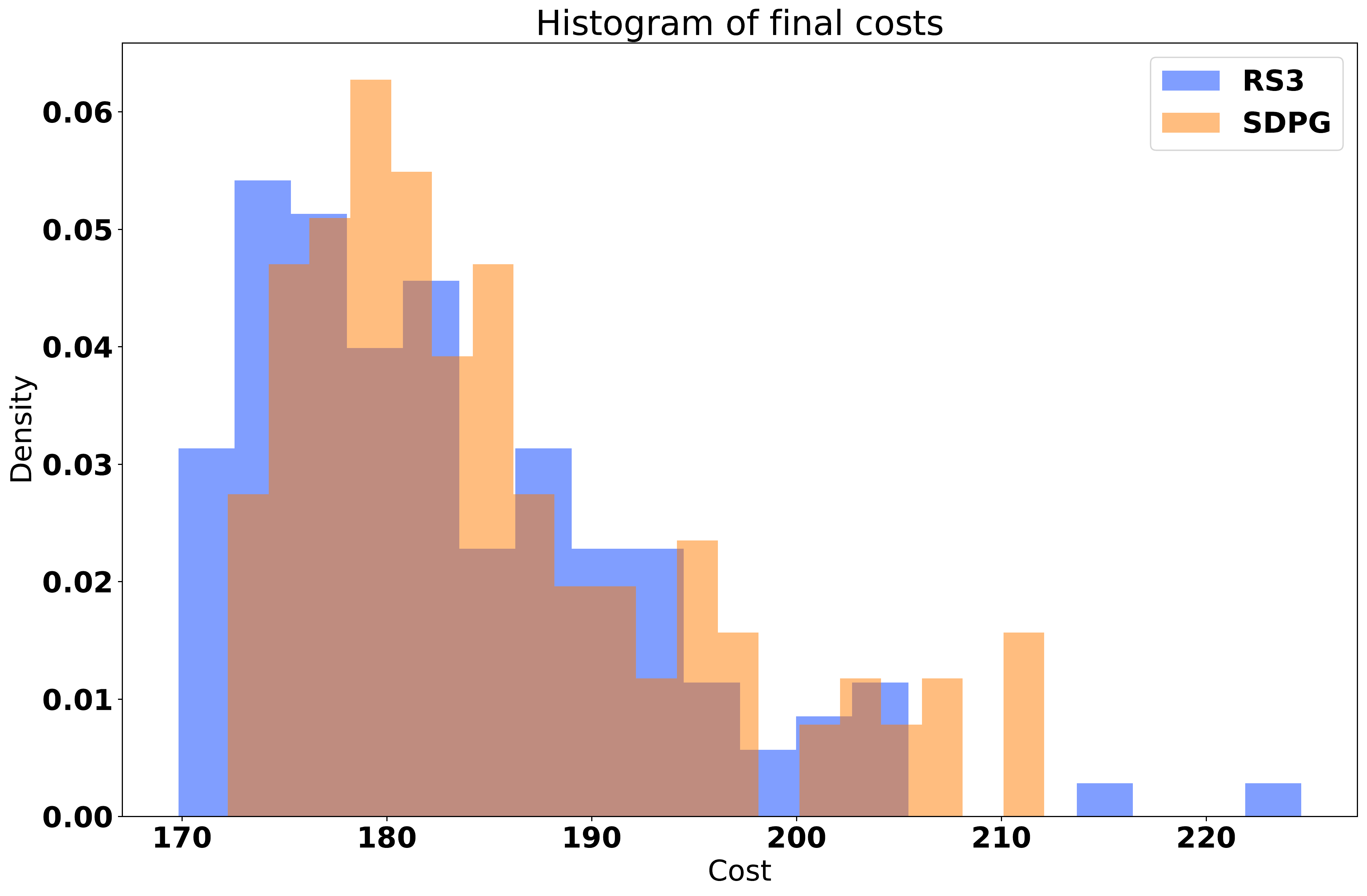}
    \includegraphics[width=0.49\textwidth]{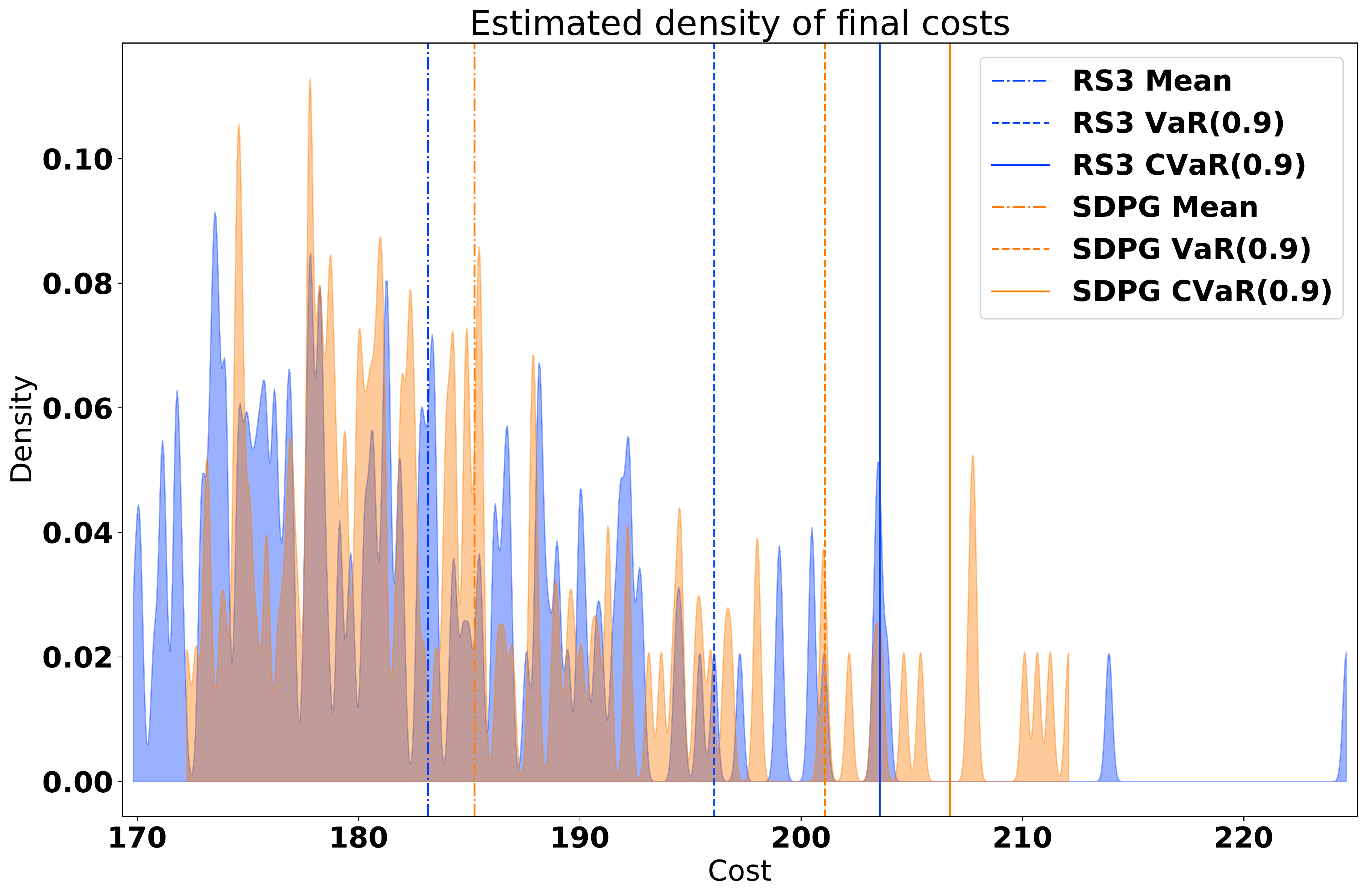}
    \caption{\textit{Top}: Control and state trajectories from the Pendulum problem with noise in the control channel. The control noise is $\mathcal{N}(0, 3.0^2).$ \textit{Bottom}: Cost histogram and the estimated p.d.f. of it. }
    \label{fig:comparison/SDPG/Pendulum/3.0}
\end{figure}

Four different levels of control noise are tested in the pendulum simulations comparing \ac{RS3} against \ac{SDPG}: $\mathcal{N}(0, 0.3^2), \mathcal{N}(0, 1.0^2), \mathcal{N}(0, 2.0^2),$ and $\mathcal{N}(0, 3.0^2)$.
The mass of the pendulum is set to $\SI{1}{\kg}$ and its length is set to $\SI{1.0}{\m}$.
The controls are constrained by the box constraints $\abs{u_t} < 10$.
The results can be found in the figures below.
As discussed in the main paper and can be observed from \Cref{fig:comparison/SDPG/Pendulum/0.3} to \Cref{fig:comparison/SDPG/Pendulum/3.0}, \ac{RS3} outperforms \ac{SDPG} in terms of the Mean, VaR, and CVaR of the final costs in all cases.

\newpage

\subsubsection{Cartpole}
\begin{figure}[H]
    \centering
    \includegraphics[width=0.49\textwidth]{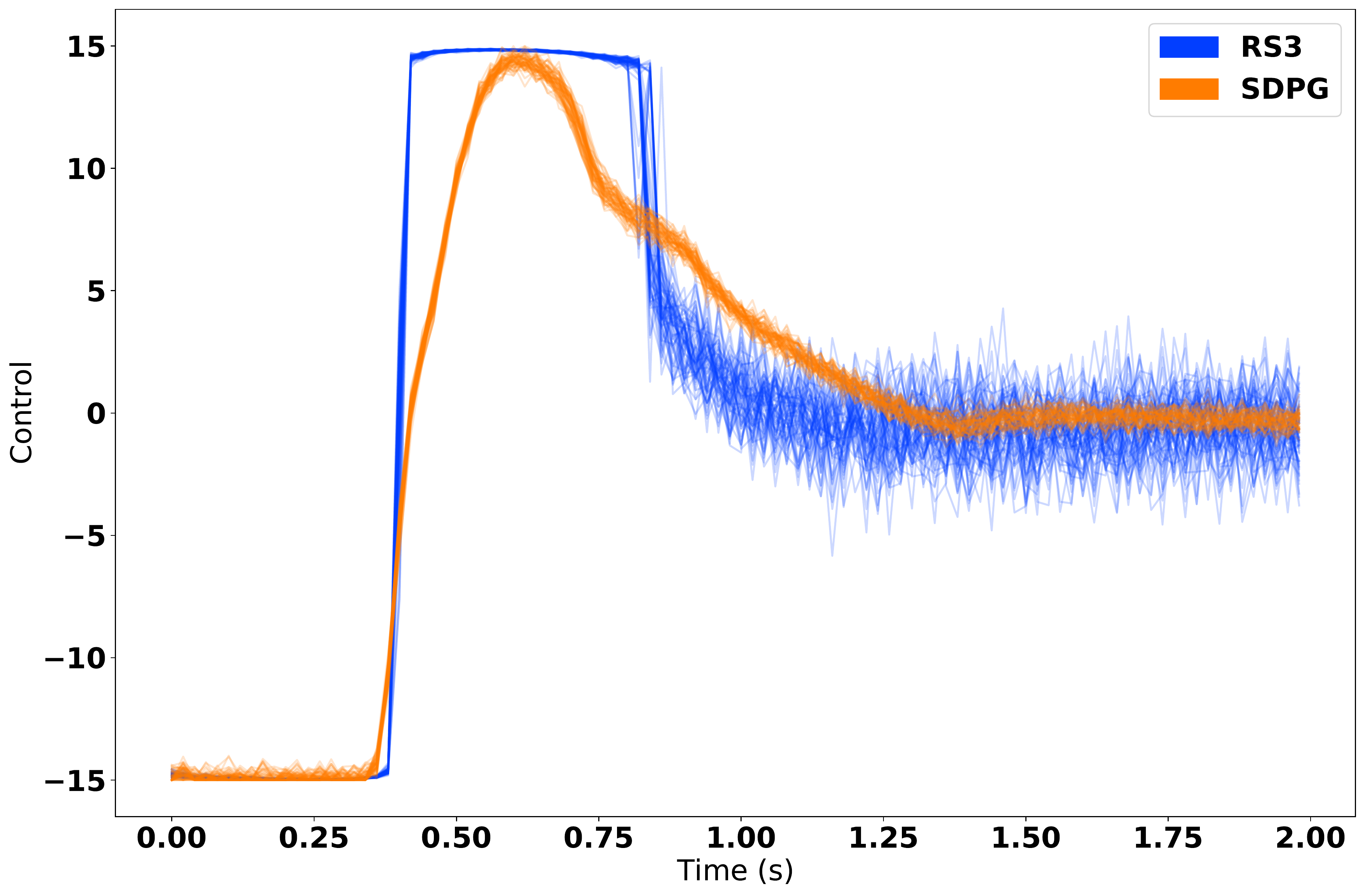}
    \includegraphics[width=0.49\textwidth]{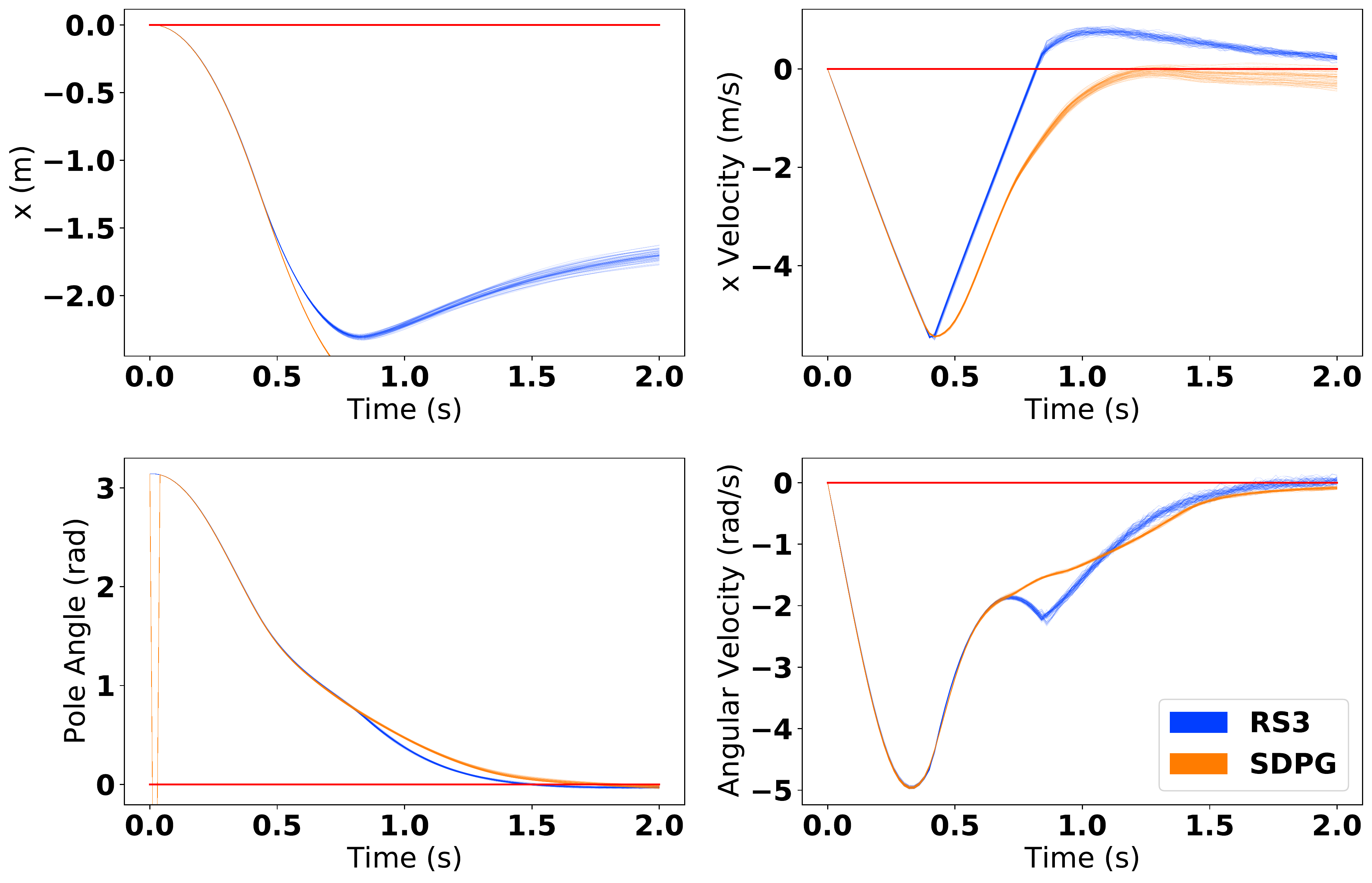}
    \includegraphics[width=0.49\textwidth]{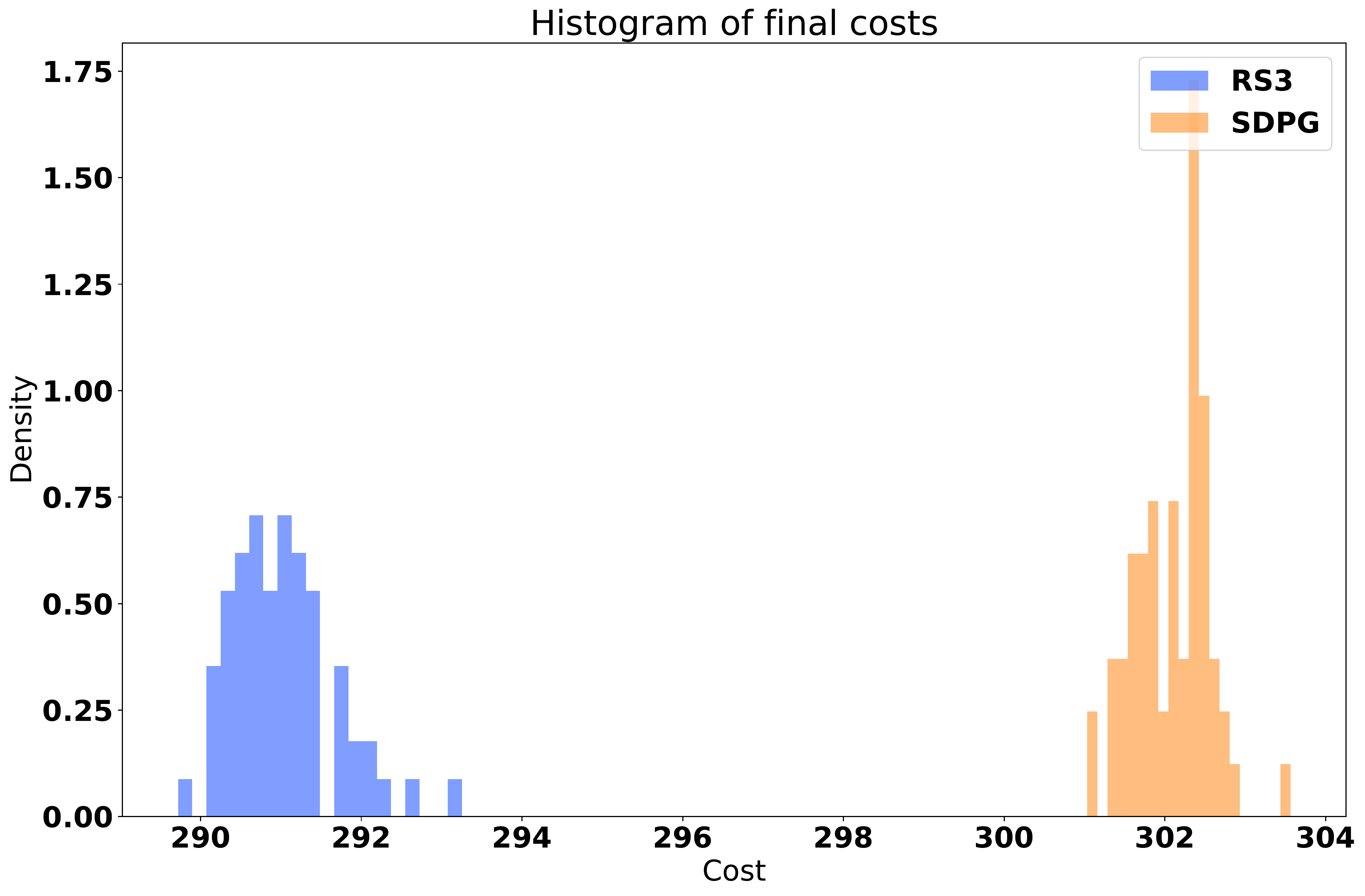}
    \includegraphics[width=0.49\textwidth]{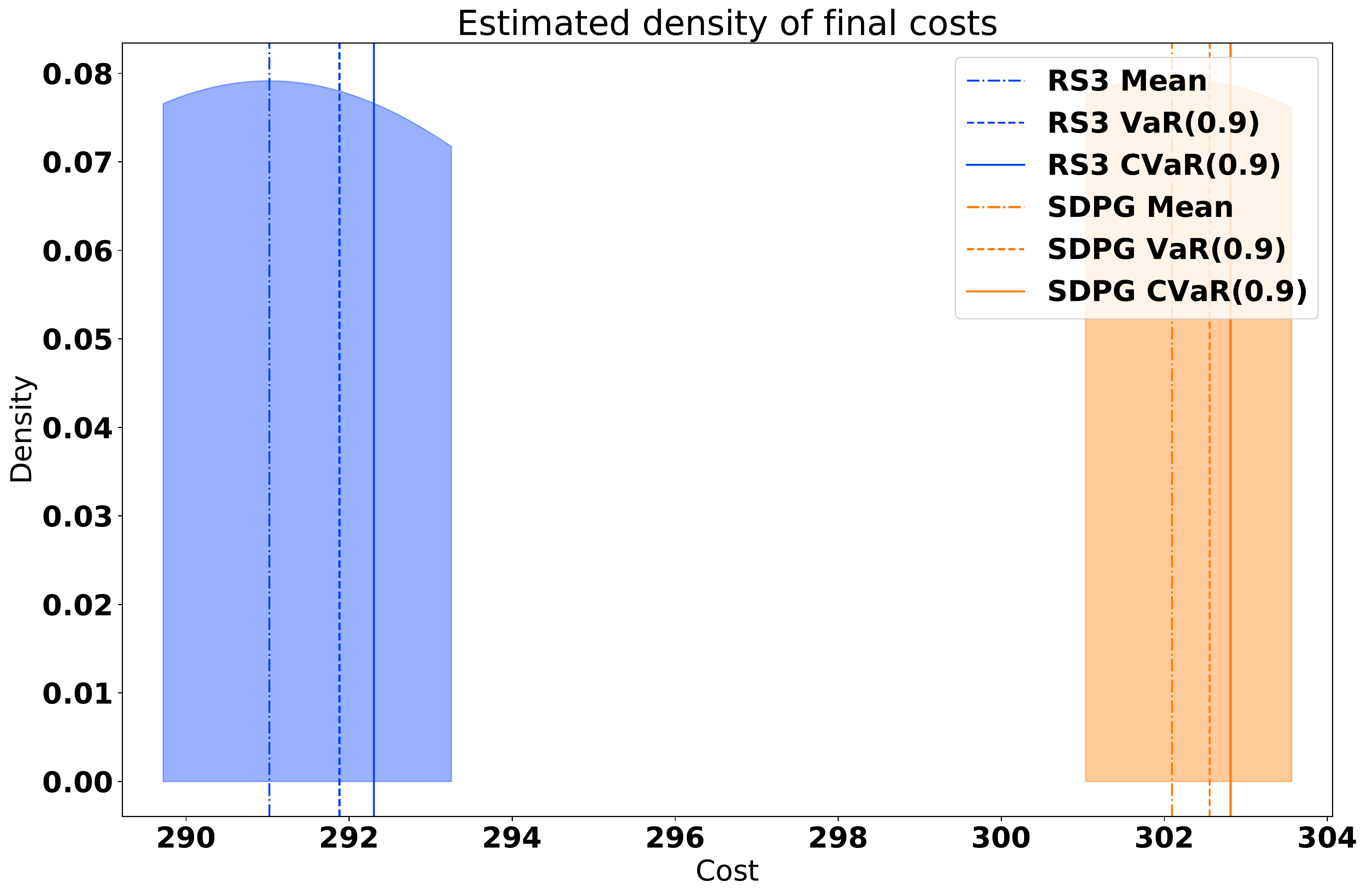}
    \caption{\textit{Top}: Control and state trajectories from the Cartpole problem with noise in the control channel. The control noise is $\mathcal{N}(0, 0.3^2).$ \textit{Bottom}: Cost histogram and the estimated p.d.f. of it. }
    \label{fig:comparison/SDPG/Cartpole/0.3}
\end{figure}

\begin{figure}[H]
    \centering
    \includegraphics[width=0.49\textwidth]{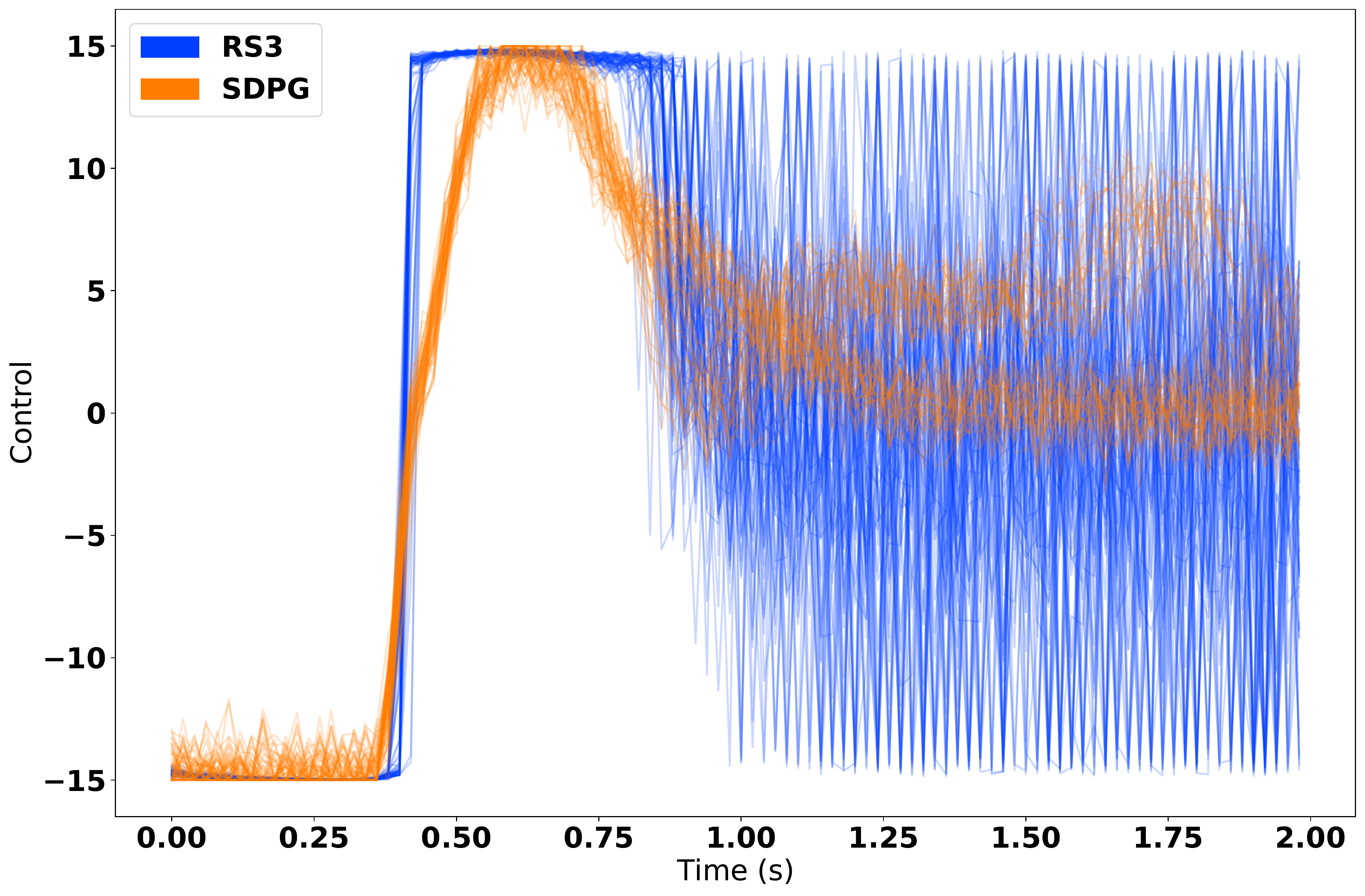}
    \includegraphics[width=0.49\textwidth]{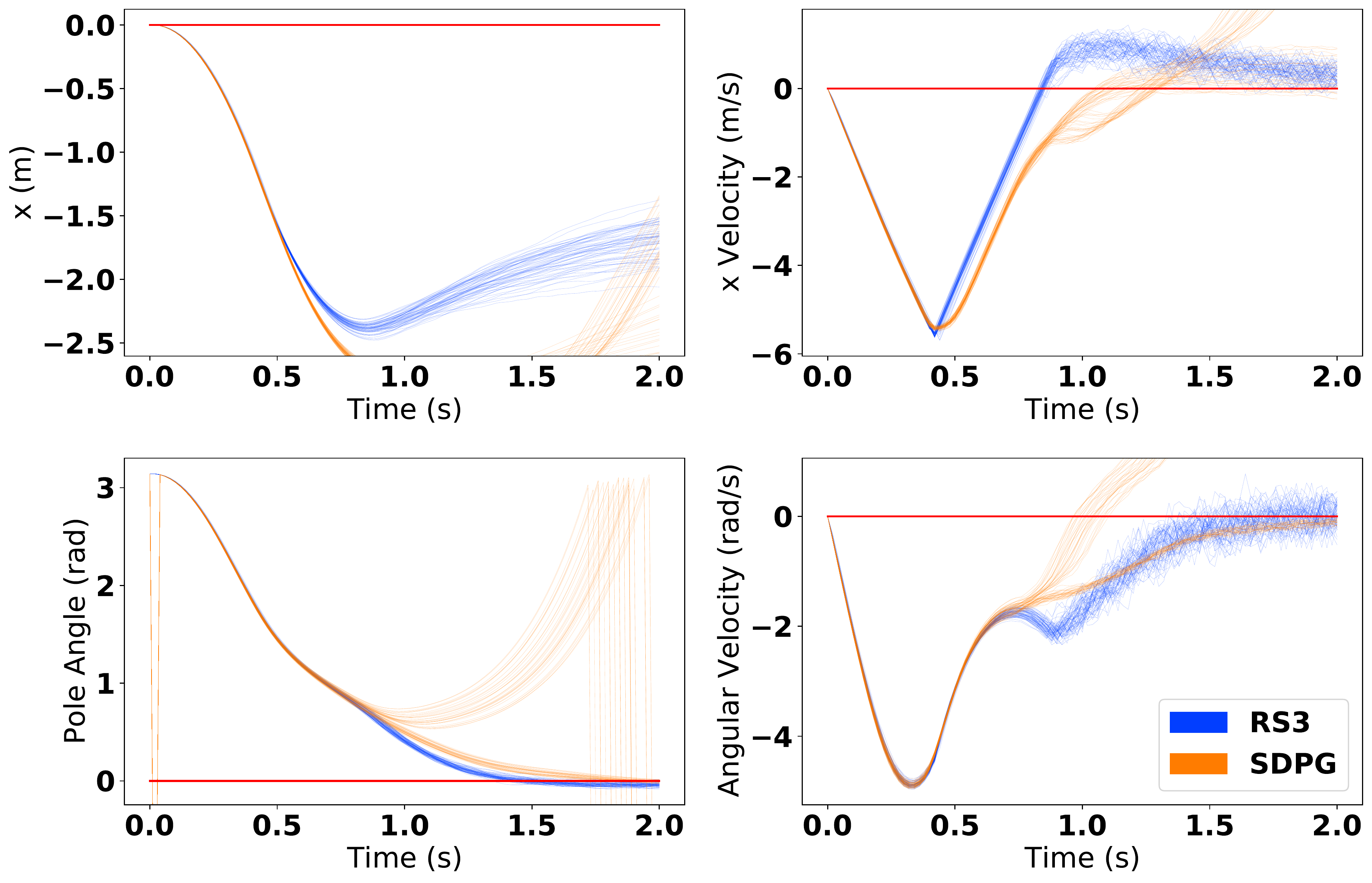}
    \includegraphics[width=0.49\textwidth]{figures/comparison/SDPG/Cartpole/1.0/cost_histogram.pdf}
    \includegraphics[width=0.49\textwidth]{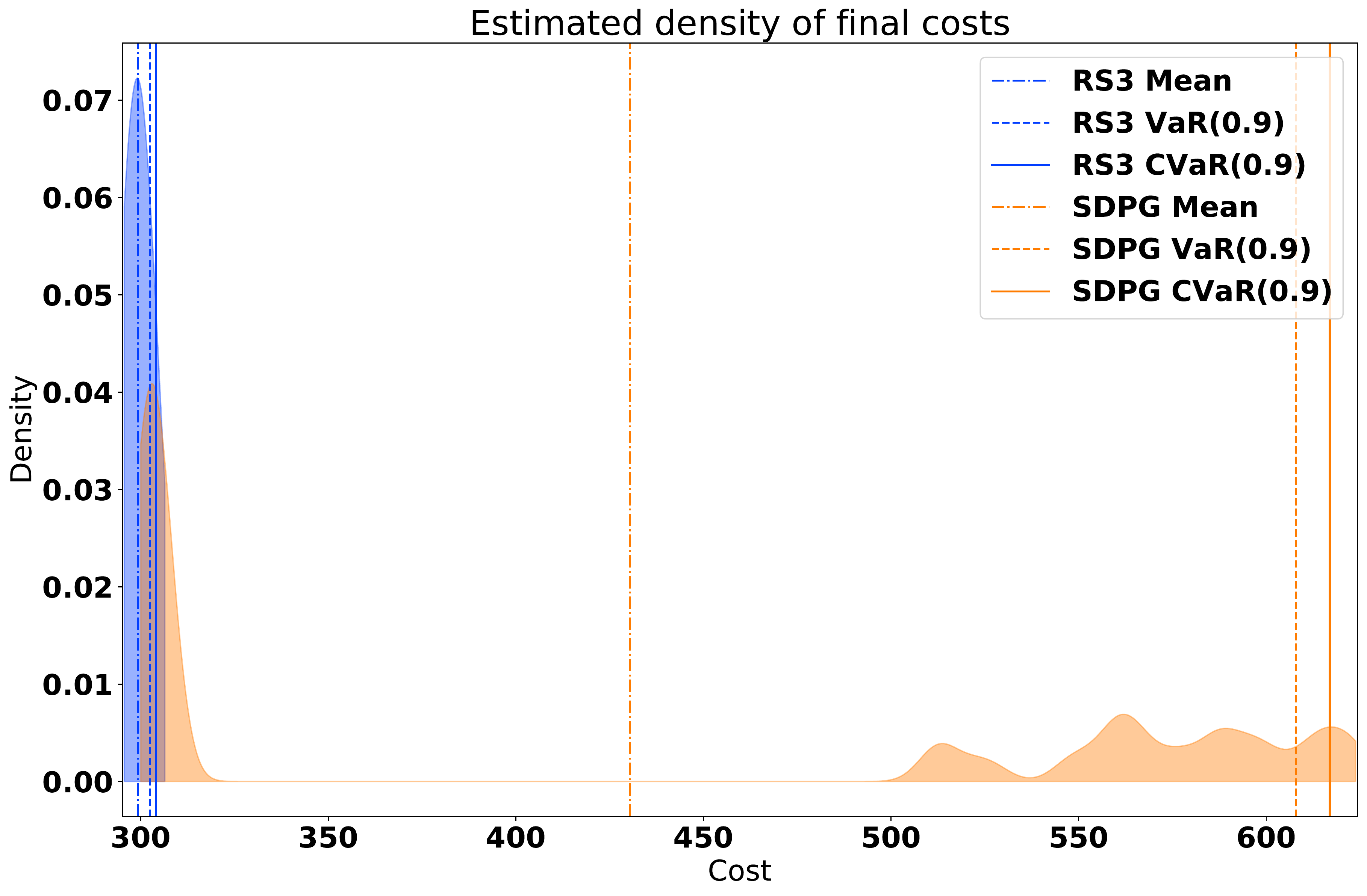}
    \caption{\textit{Top}: Control and state trajectories from the Cartpole problem with noise in the control channel. The control noise is $\mathcal{N}(0, 1.0^2).$ \textit{Bottom}: Cost histogram and the estimated p.d.f. of it. }
    \label{fig:comparison/SDPG/Cartpole/1.0}
\end{figure}

\begin{figure}[H]
    \centering
    \includegraphics[width=0.49\textwidth]{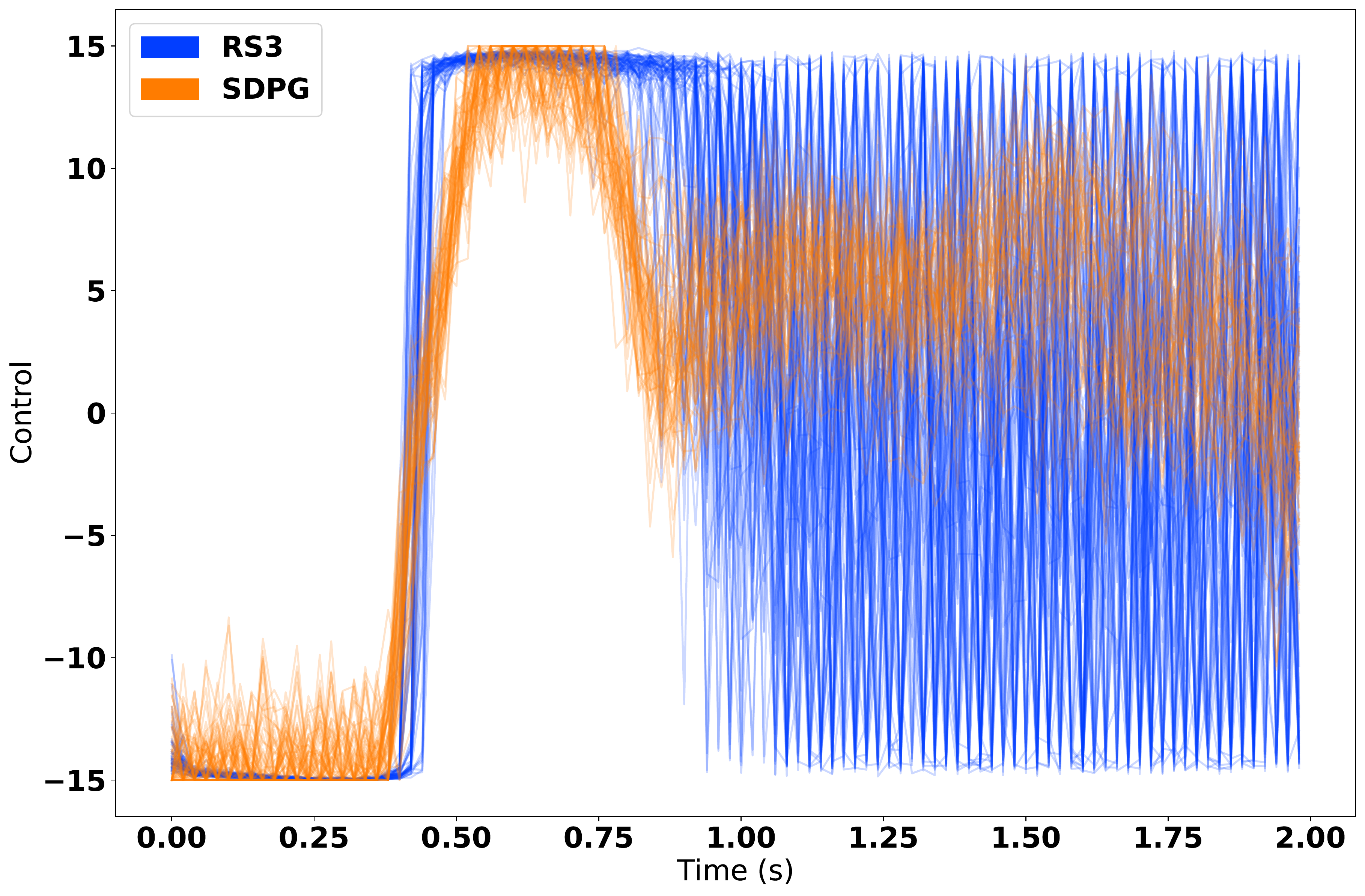}
    \includegraphics[width=0.49\textwidth]{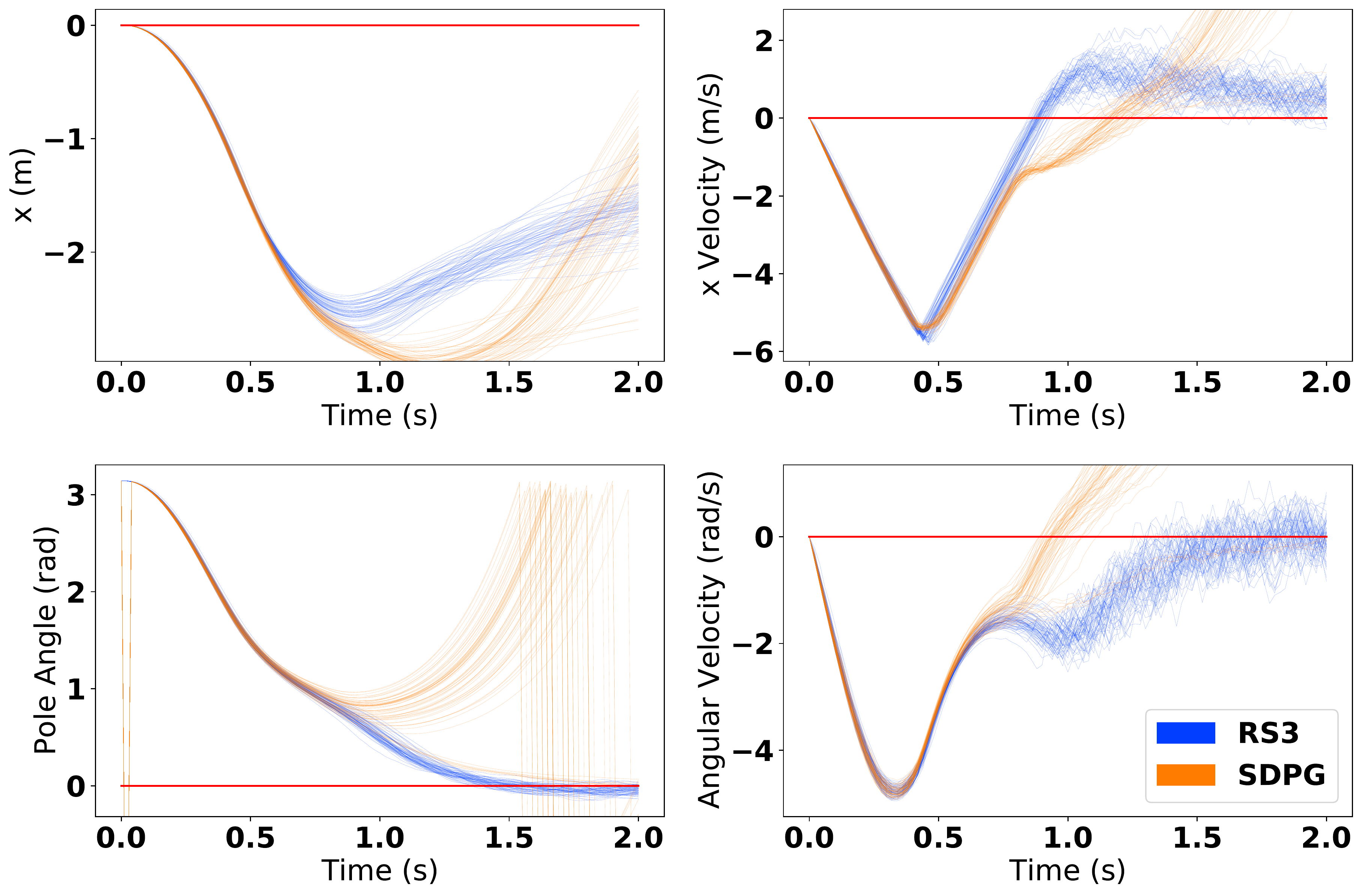}
    \includegraphics[width=0.49\textwidth]{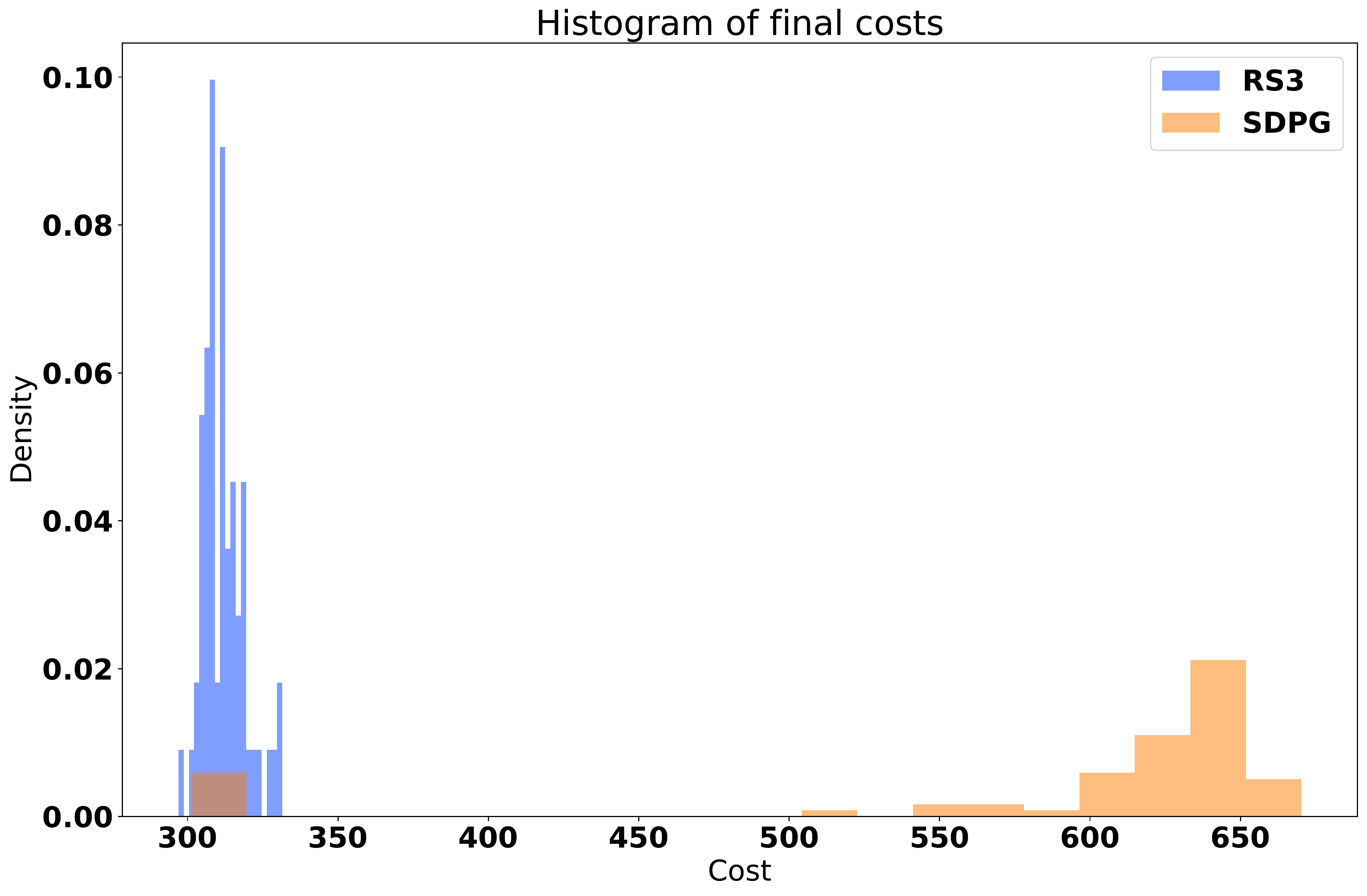}
    \includegraphics[width=0.49\textwidth]{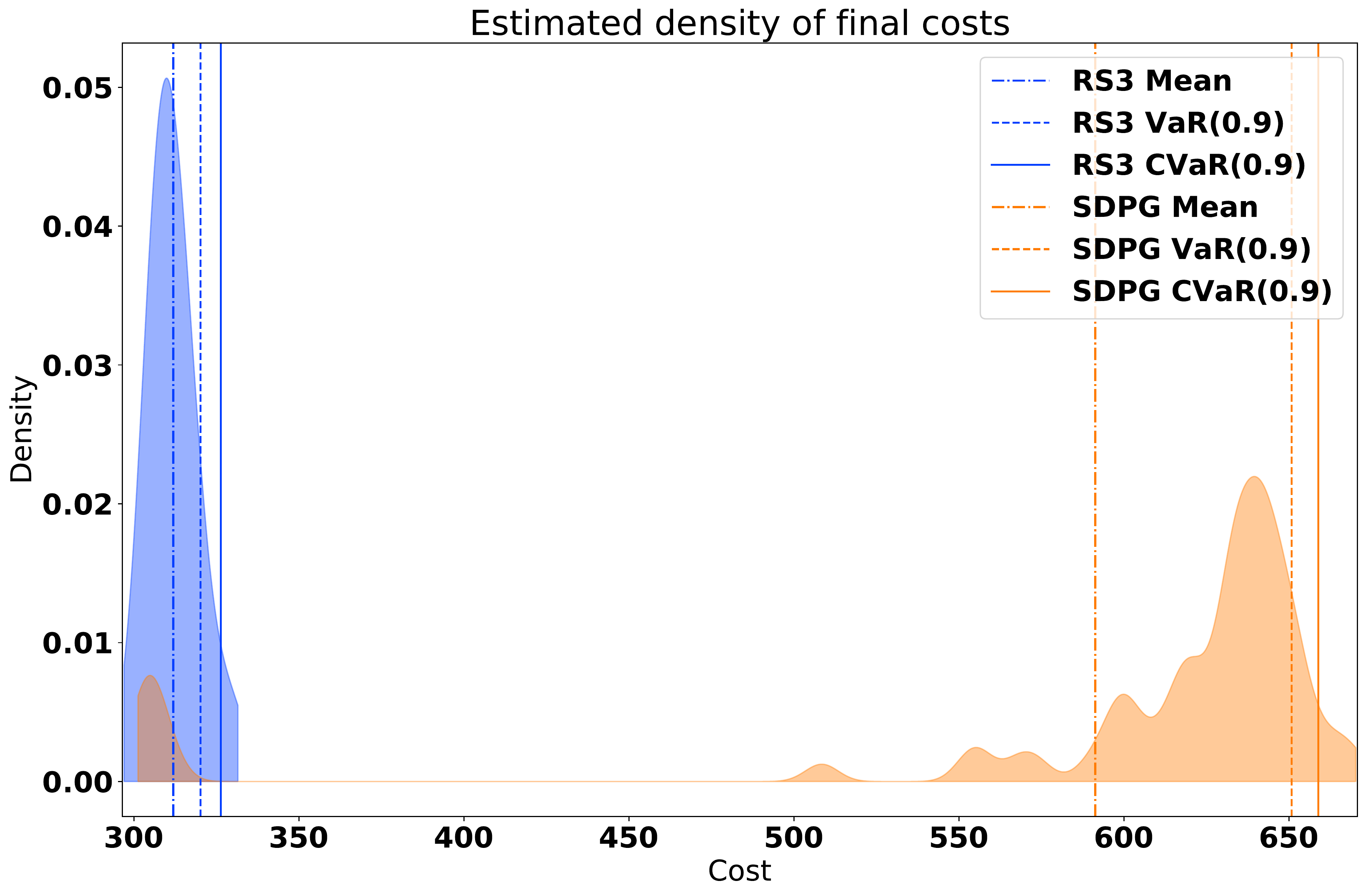}
    \caption{\textit{Top}: Control and state trajectories from the Cartpole problem with noise in the control channel. The control noise is $\mathcal{N}(0, 2.0^2).$ \textit{Bottom}: Cost histogram and the estimated p.d.f. of it. }
    \label{fig:comparison/SDPG/Cartpole/2.0}
\end{figure}

Similarly, 3 different levels of additive control noise were used to compare \ac{RS3} and \ac{SDPG} in the cartpole problem: $\mathcal{N}(0, 0.3^2), \mathcal{N}(0, 1.0^2),$ and $\mathcal{N}(0, 2.0^2)$.
The mass of the cart are set to $\SI{1}{\kg}$, the mass of the pole $\SI{0.1}{\kg}$ and its length $\SI{0.5}{\m}$.
The controls are constrained by the box constraints $\abs{u_t} < 15$.
The results can be found in \Cref{fig:comparison/SDPG/Cartpole/0.3} to \Cref{fig:comparison/SDPG/Cartpole/2.0}.

Although the control becomes spiky with larger noise, \ac{RS3} is able to accomplish the task, whereas \ac{SDPG} fails the task most of the time, as seen in \Cref{fig:comparison/SDPG/Cartpole/1.0} and \Cref{fig:comparison/SDPG/Cartpole/2.0}.
These results show that \ac{RS3} does a better job of risk-sensitive control on higher order systems compared to \ac{SDPG}.

\newpage

\section{Belief Space Risk Sensitive Control}
\subsection{Pendulum}
We test \ac{RS3} on Pendulum in the case of uncertainty in the initial conditions and in the case of stochastic dynamics where the true state of the system is fully observable.
The mass of the pendulum is set to $\SI{1}{\kg}$ and its length is set to $\SI{1.0}{\m}$.
The controls are constrained by the box constraints $\abs{u_t} < 10$.
Measurements of the state $x$ are corrupted with additive noise $\xi \sim \mathcal{N}(0, \textrm{diag}([0.7, 0.3]))$.
In the case of uncertain initial states, the initial states $x_0 = [\theta, \dot{\theta} ]$ are drawn from a normal distribution with mean $[\pi, 0]$ and covariance matrix $\textrm{diag}([0.5, 0.5])$.
In the case of stochastic dynamics, the controls are corrupted by noise with distribution $\mathcal{N}(0, 3.0^2)$.

\begin{figure}[htb!]
    \centering
    \includegraphics[width=\textwidth]{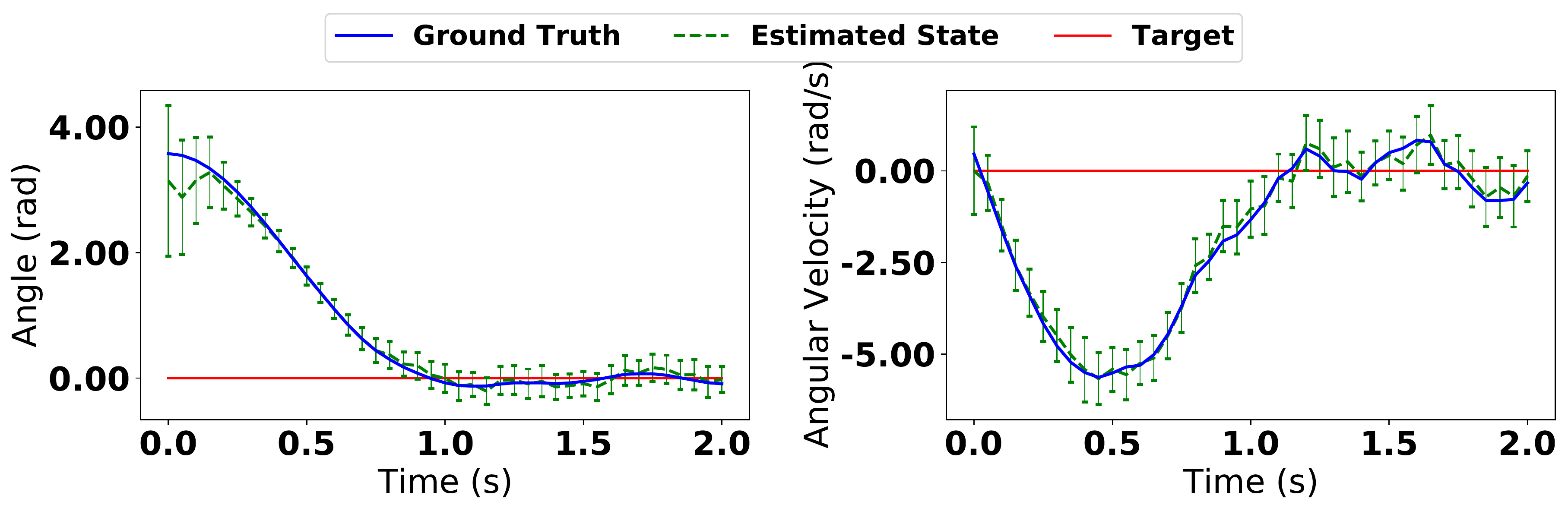}
    \includegraphics[width=\textwidth]{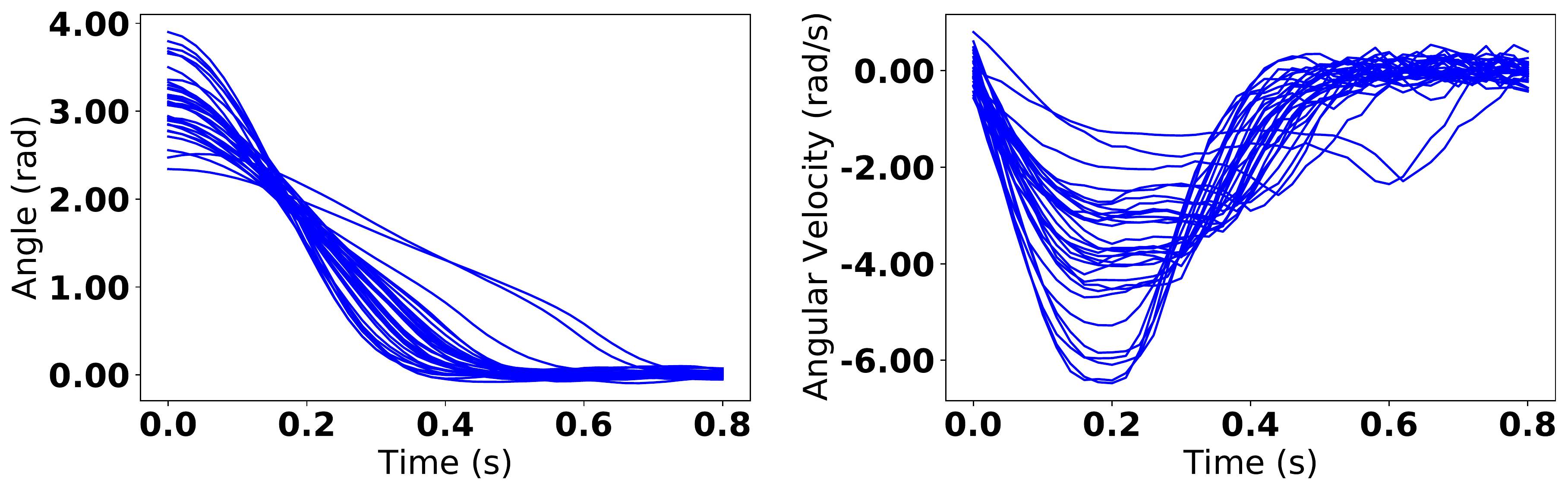}
    \caption{Trajectories from the Pendulum problem with is uncertainty in the initial conditions.
        \textit{Top}: A comparison of the posterior state from the particle filter and the ground truth from a single trajectory.
        \textit{Bottom}: Different trajectory realizations of the initial condition.}
    \label{fig:pf_pendulum/est_state_comp}
\end{figure}

\begin{figure}[htb!]
    \centering
    \includegraphics[width=\textwidth]{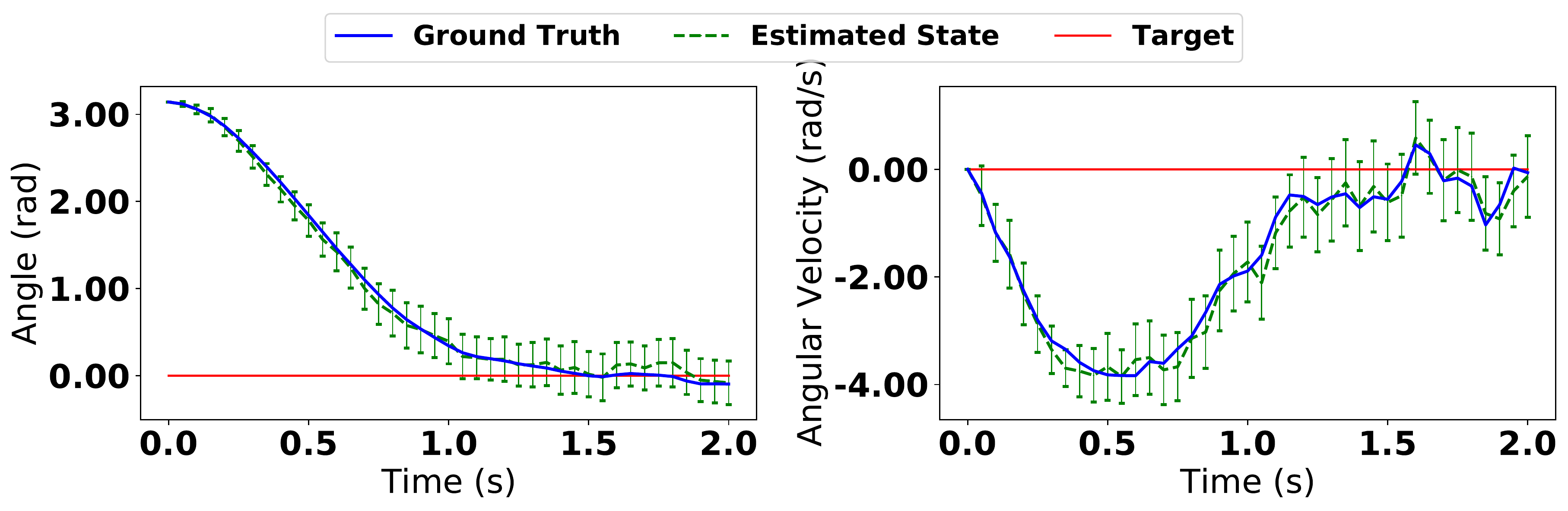}
    \includegraphics[width=\textwidth]{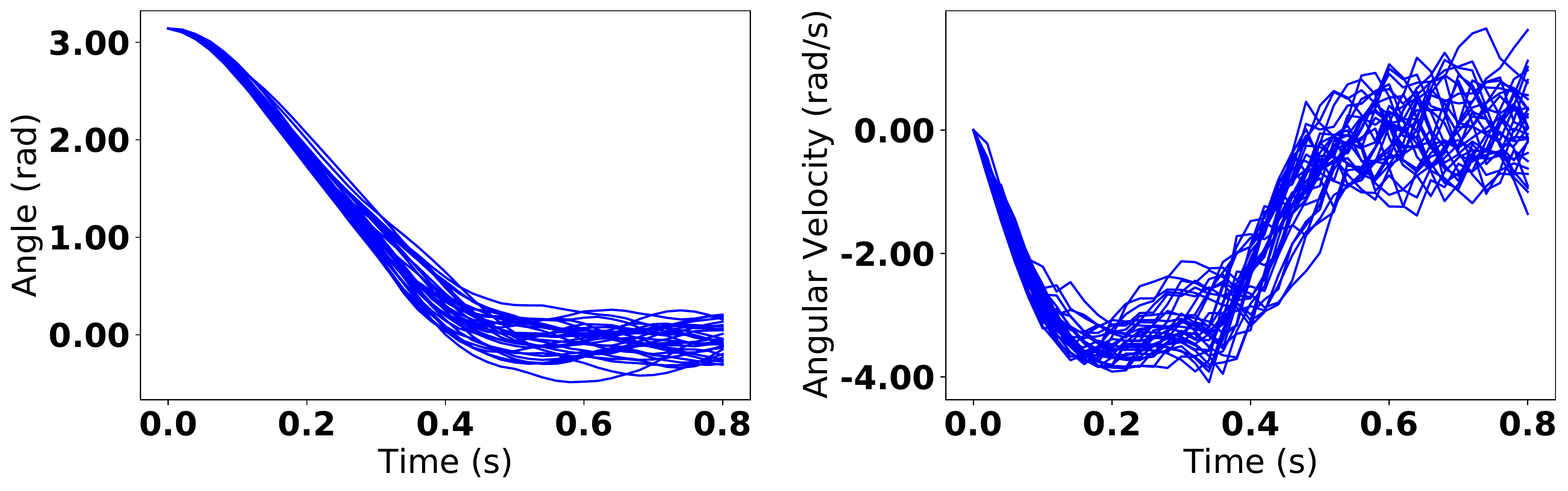}
    \caption{Trajectories from the Pendulum problem with stochastic dynamics.
        \textit{Top}: A comparison of the posterior state from the particle filter and the ground truth from a single trajectory.
        \textit{Bottom}: Different trajectory realizations of the stochastic dynamics.}
    \label{fig:pf_control_noise/pendulum}
\end{figure}

The results for uncertain initial conditions and stochastic dynamics can be found in \Cref{fig:pf_pendulum/est_state_comp} and \Cref{fig:pf_control_noise/pendulum} respectively.

\newpage

\subsection{Cartpole}
Similarly, we test \ac{RS3} on Cartpole in the case of uncertainty in the initial conditions and in the case of stochastic dynamics.
The mass of the cart are set to $\SI{1}{\kg}$, the mass of the pole $\SI{0.1}{\kg}$ and its length $\SI{0.5}{\m}$.
The controls are constrained by the box constraints $\abs{u_t} < 15$.
Measurements of the state $x$ are corrupted with additive noise $\xi \sim \mathcal{N}(0, \textrm{diag}([1, 1, 0.25, 0.25]))$.

The initial states $x_0 = [x, \dot{x}, \theta, \dot{\theta} ]$ are drawn from a normal distribution with mean $[0, 0, \pi, 0]$ and covariance matrix $\textrm{diag}([0.5, 0.5, 0.08, 0.05])$

\begin{figure}[htb!]
    \centering
    \includegraphics[width=\textwidth]{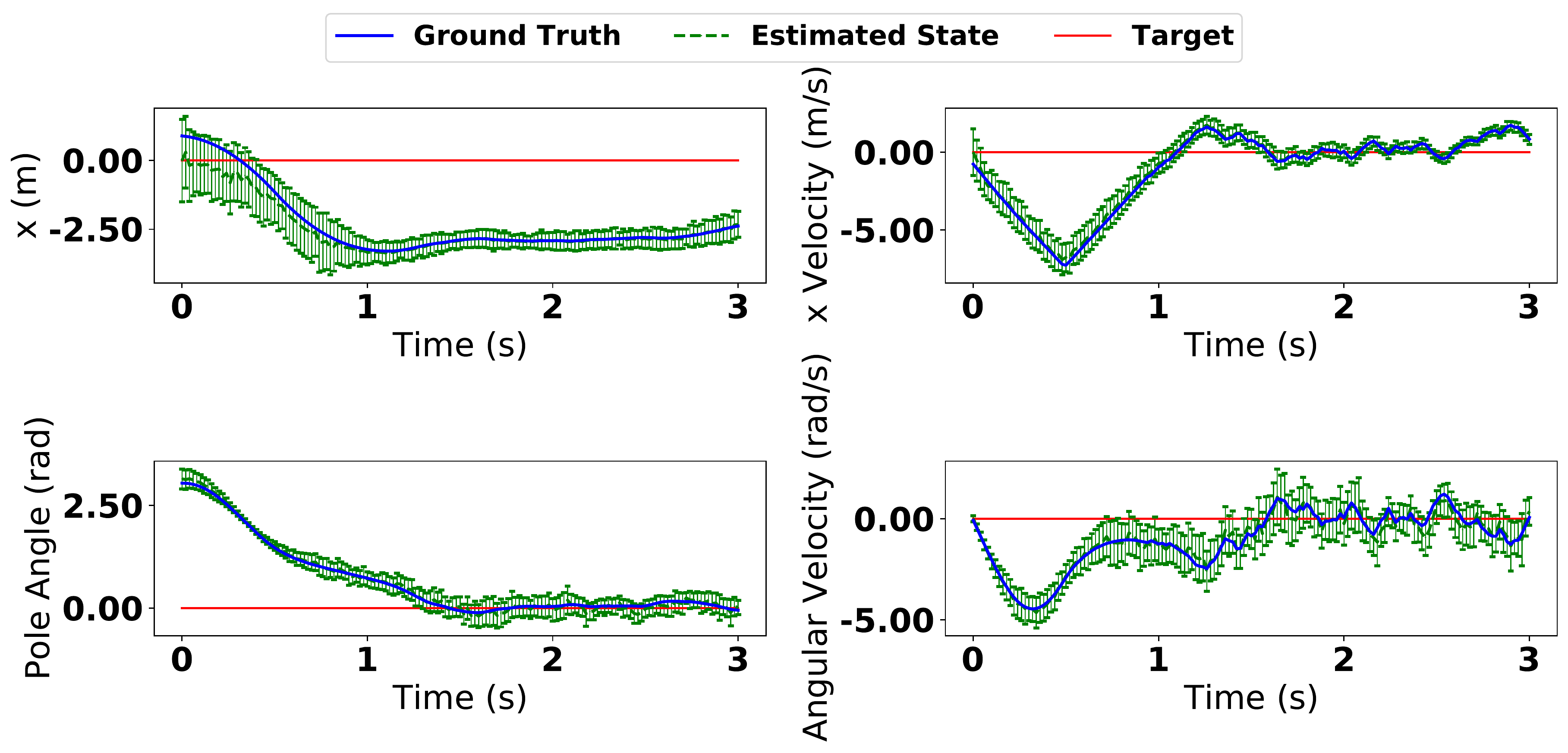}
    \includegraphics[width=\textwidth]{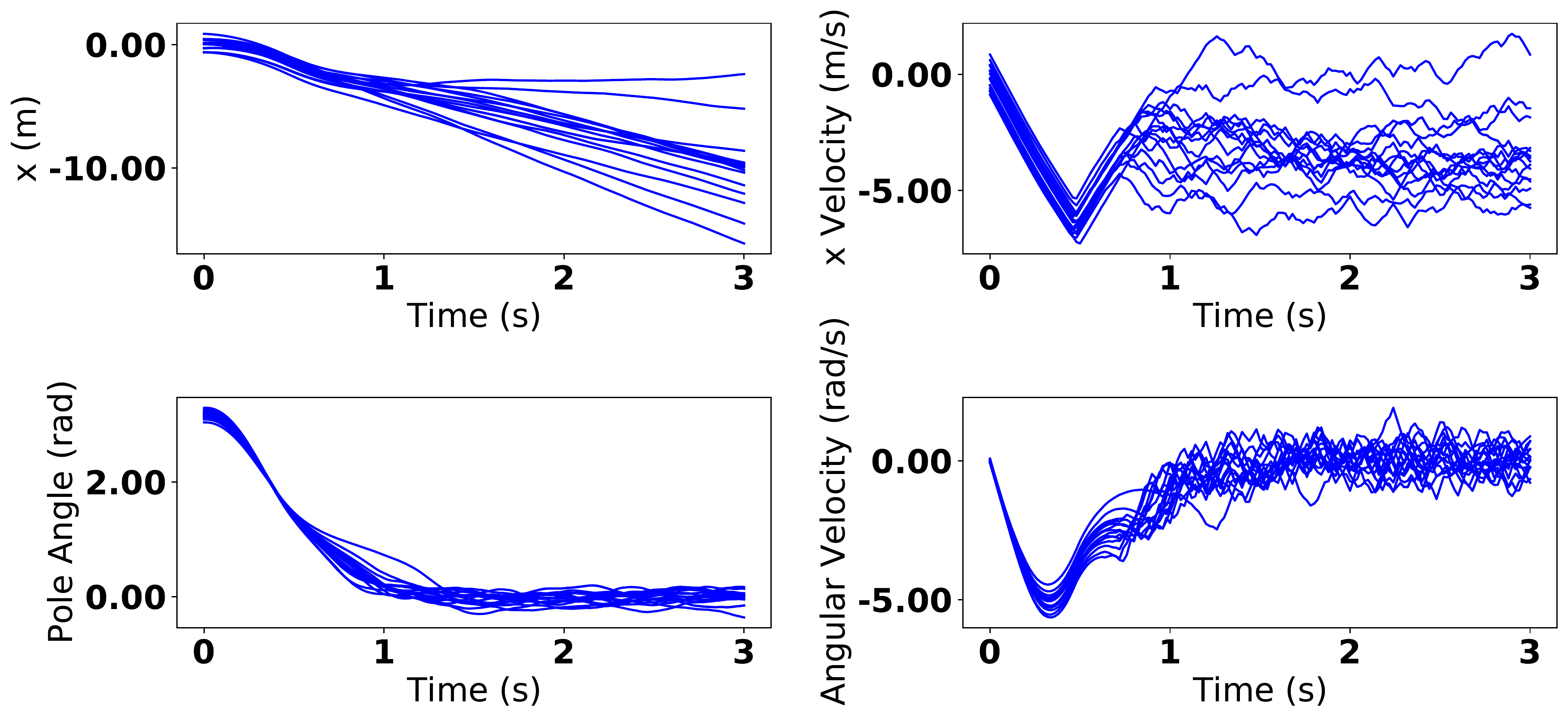}
    \caption{Trajectories from the Cartpole problem with is uncertainty in the initial conditions.
        \textit{Top}: A comparison of the posterior state from the particle filter and the ground truth from a single trajectory.
        \textit{Bottom}: Different trajectory realizations of the initial condition.}
    \label{fig:pf_cartpole/est_state_comp}
\end{figure}

\begin{figure}[htb!]
    \centering
    \includegraphics[width=\textwidth]{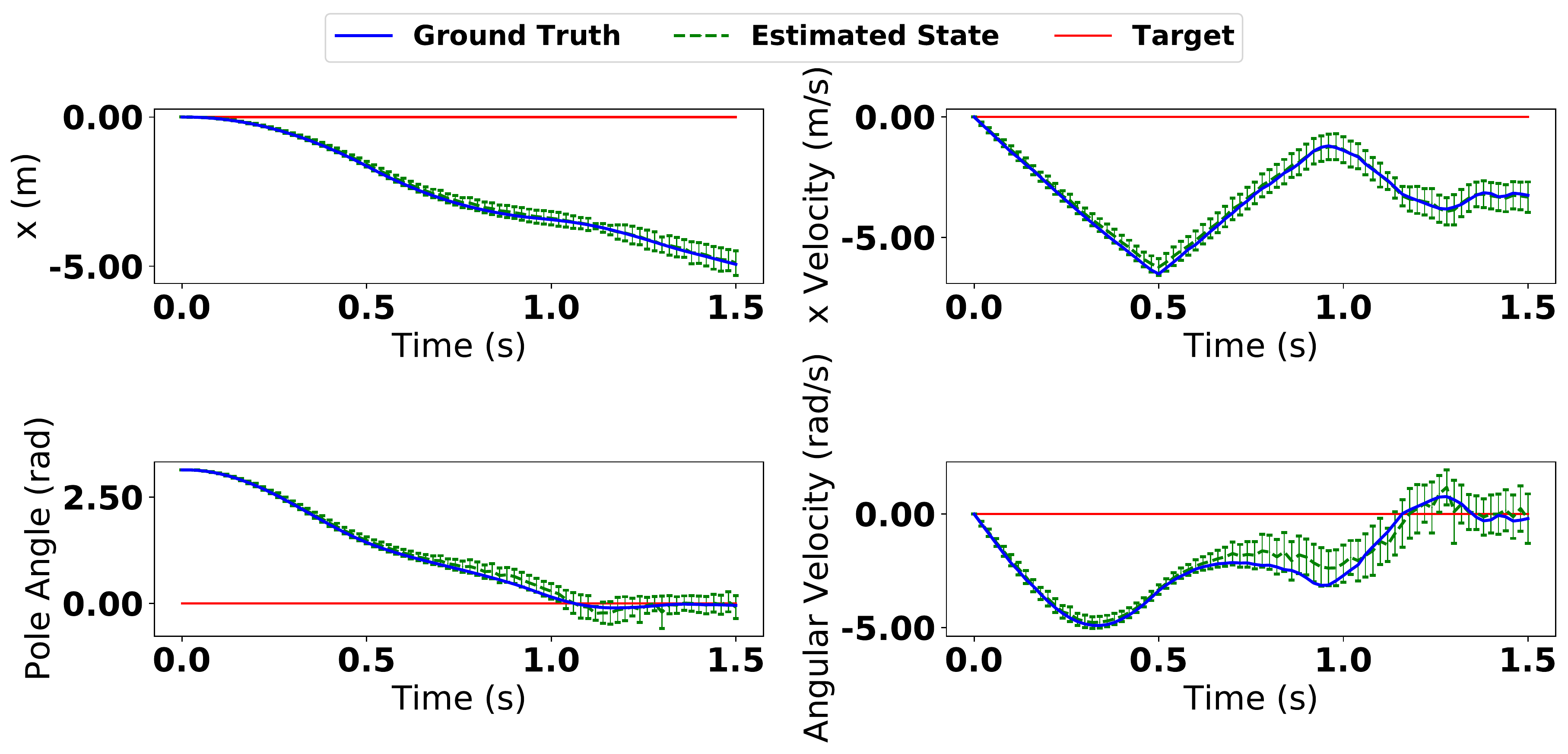}
    \includegraphics[width=\textwidth]{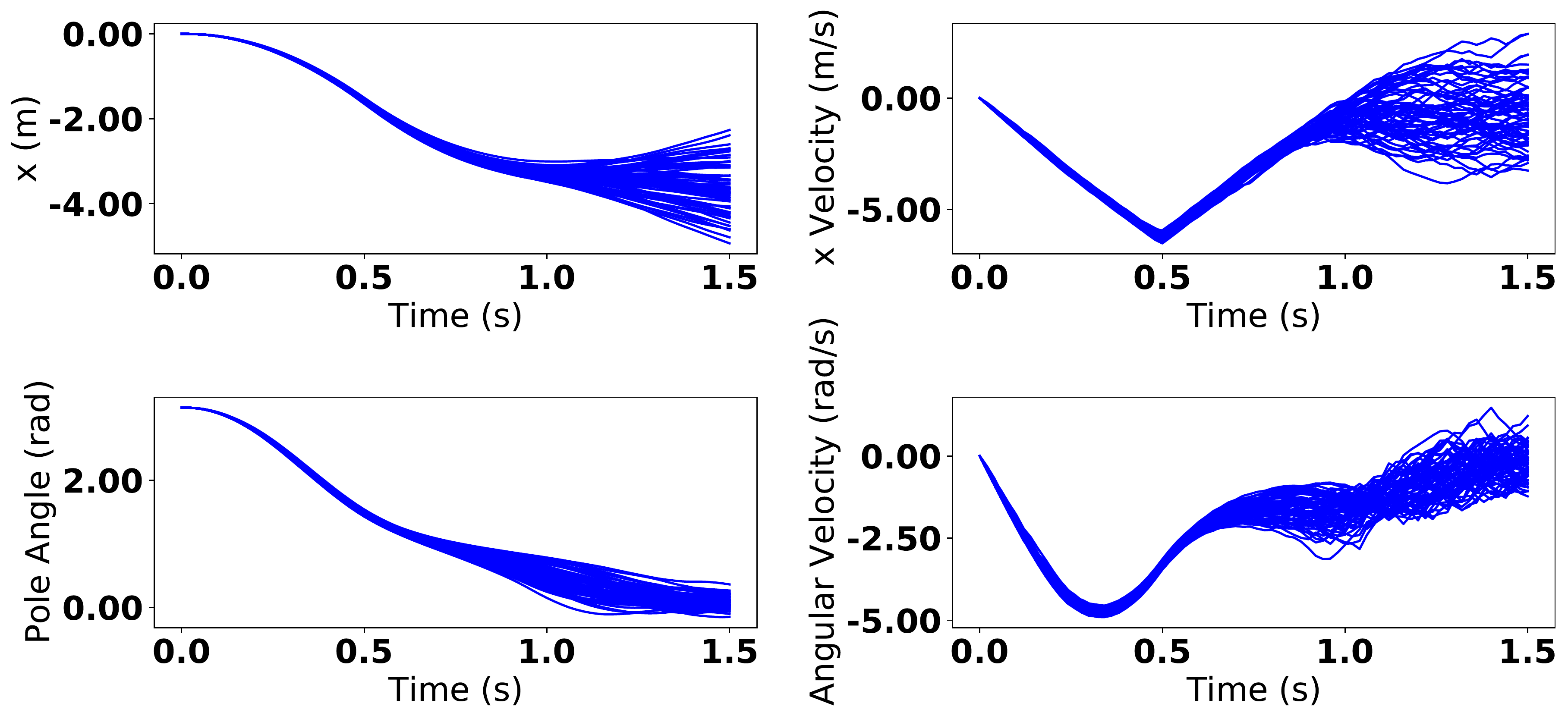}
    \caption{Trajectories from the Cartpole problem with stochastic dynamics.
        \textit{Top}: A comparison of the posterior state from the particle filter and the ground truth from a single trajectory.
        \textit{Bottom}: Different trajectory realizations of the stochastic dynamics.}
    \label{fig:pf_control_noise/cartpole}
\end{figure}

The results for uncertain initial conditions and stochastic dynamics can be found in \Cref{fig:pf_cartpole/est_state_comp} and \Cref{fig:pf_control_noise/cartpole} respectively.

\newpage

\subsection{Quadcopter}
Finally, we test \ac{RS3} on the Quadcopter in the case of uncertainty in the initial conditions and stochastic dynamics.
Measurements of the state $x$ are corrupted with additive noise $\xi \sim \mathcal{N}(0, \textrm{diag}([0.1, 0.1, 0.1, 0.01, 0.01, 0.01, 0.08, 0.08, 0.08, 0.01, 0.1, 0.01]))$.
The controls are constrained by the box constraints $[0, -10, -10, -1] < u_t < [20, 10, 10, 1]$.

For the case of uncertain initial states, the states $x_0 = [x, y, z, \alpha, \beta, \gamma, \dot{x}, \dot{y}, \dot{z}, \dot{\alpha}, \dot{\beta}, \dot{\gamma}]$ are drawn from a normal distribution with zero mean and covariance matrix 
$\textrm{diag}([0.3, 0.3, 0.3, 0.2, 0.2, 0.2, 0.2, 0.2, 0.2, 0.2, 0.2, 0.2])$. 
In the case of stochastic dynamics, the controls are corrupted with noise drawn from $\mathcal{N}(0, \textrm{diag}([1, 1, 1, 0.1])$

\begin{figure}[htb!]
    \centering
    \includegraphics[width=\textwidth]{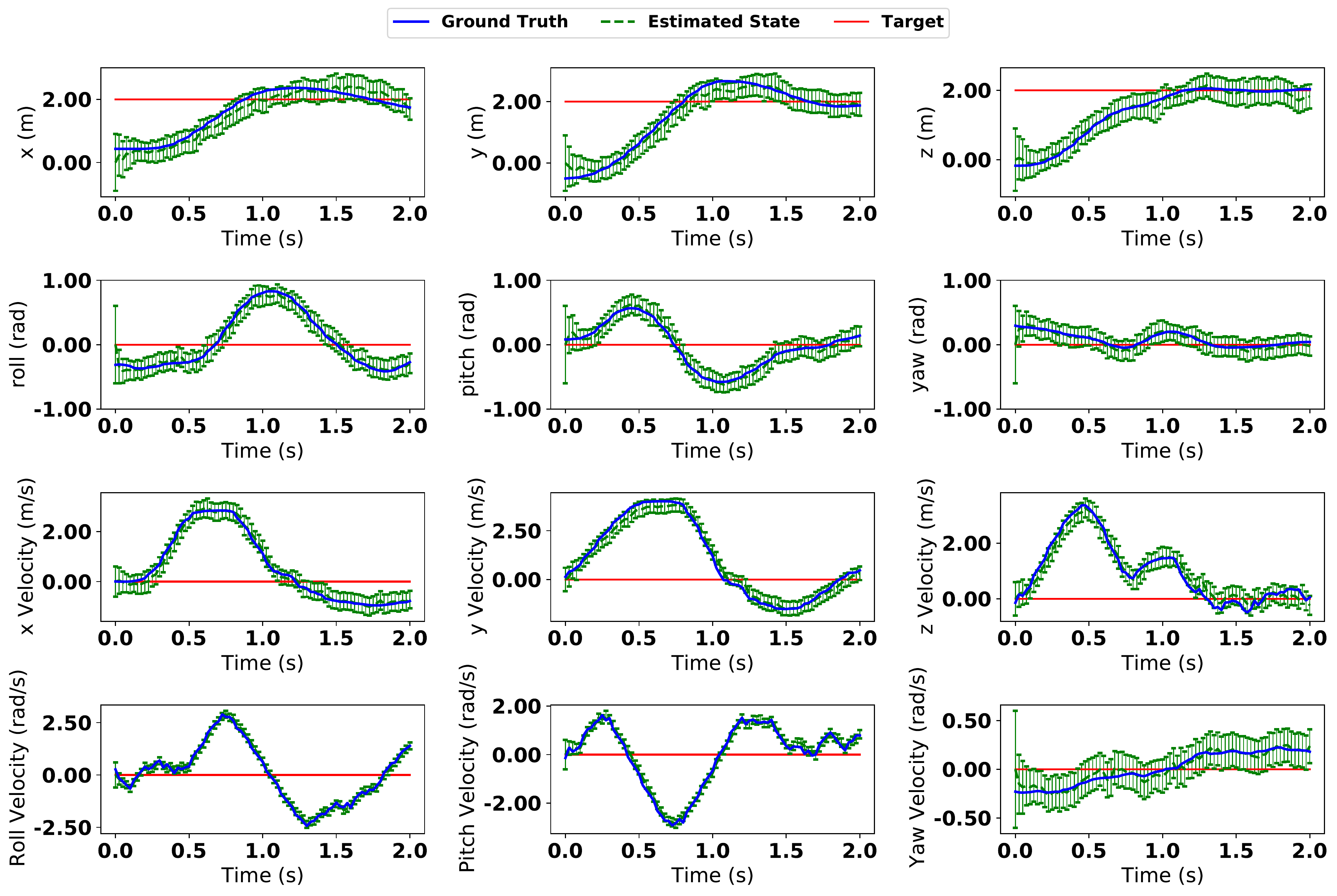}
    \includegraphics[width=\textwidth]{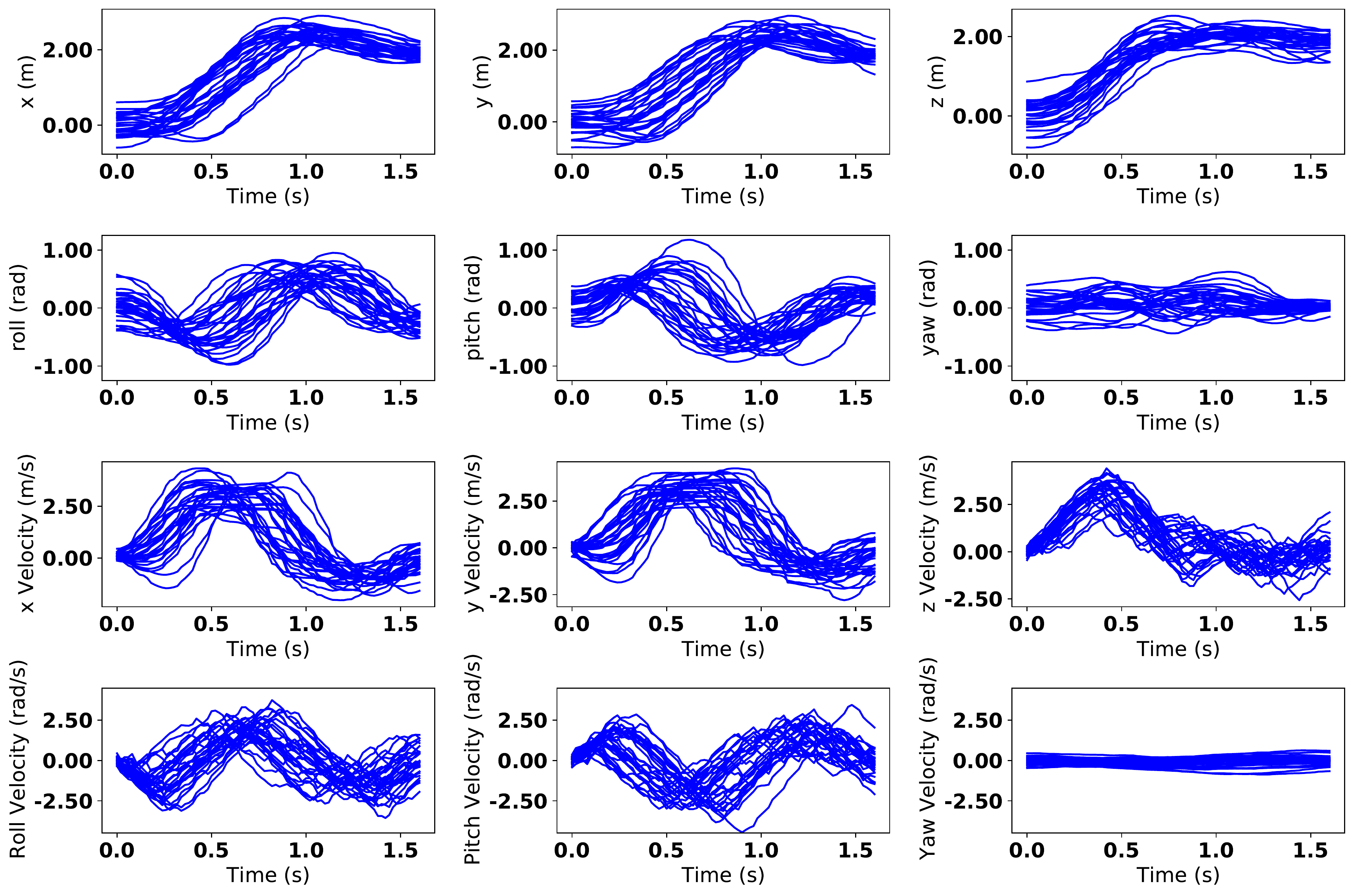}
    \caption{Trajectories from the Quadcopter problem with is uncertainty in the initial conditions.
        \textit{Top}: A comparison of the posterior state from the particle filter and the ground truth from a single trajectory.
        \textit{Bottom}: Different trajectory realizations of the initial condition.}
    \label{fig:pf_quadcopter/est_state_comp}
\end{figure}

\begin{figure}[htb!]
    \centering
    \includegraphics[width=\textwidth]{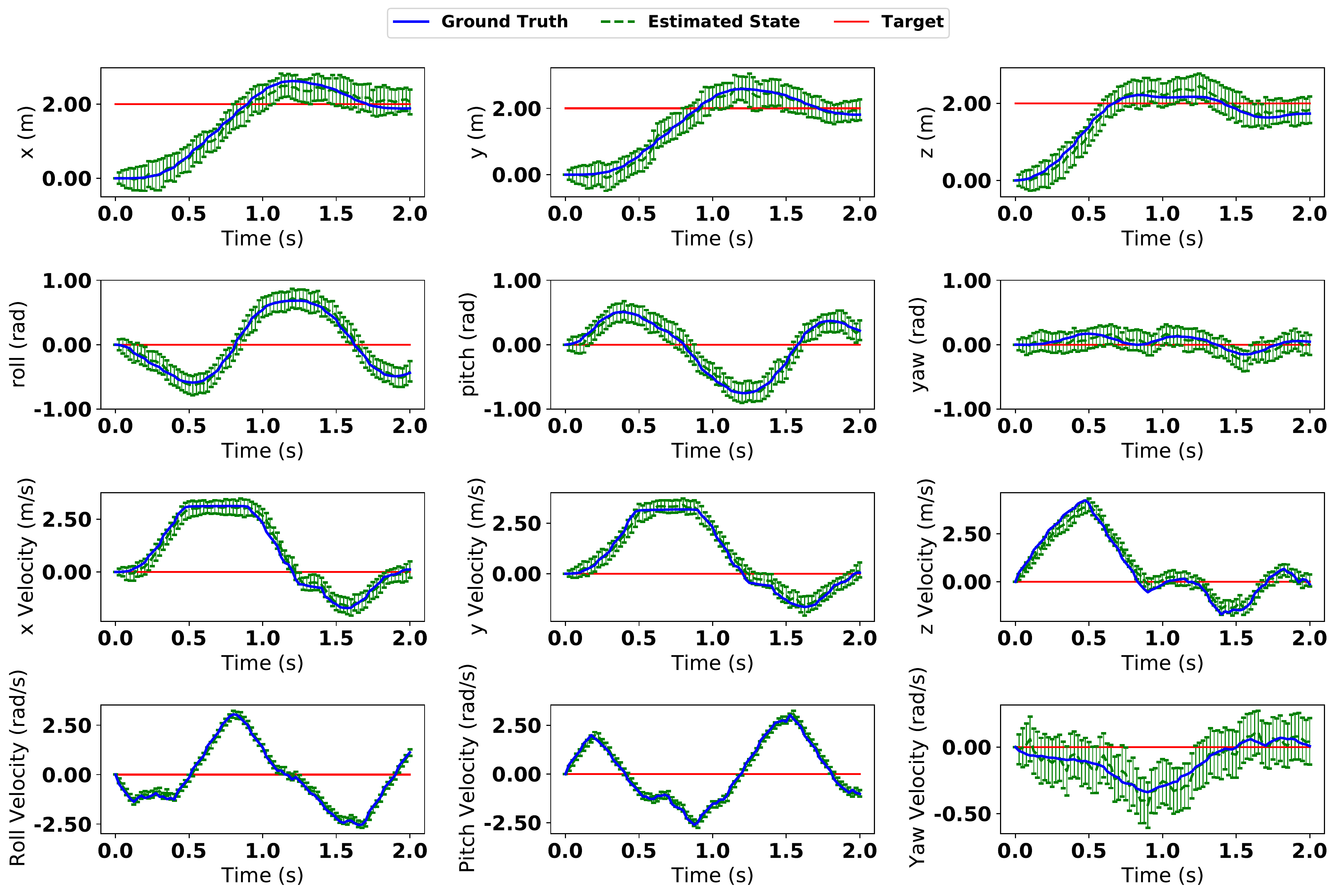}
    \includegraphics[width=\textwidth]{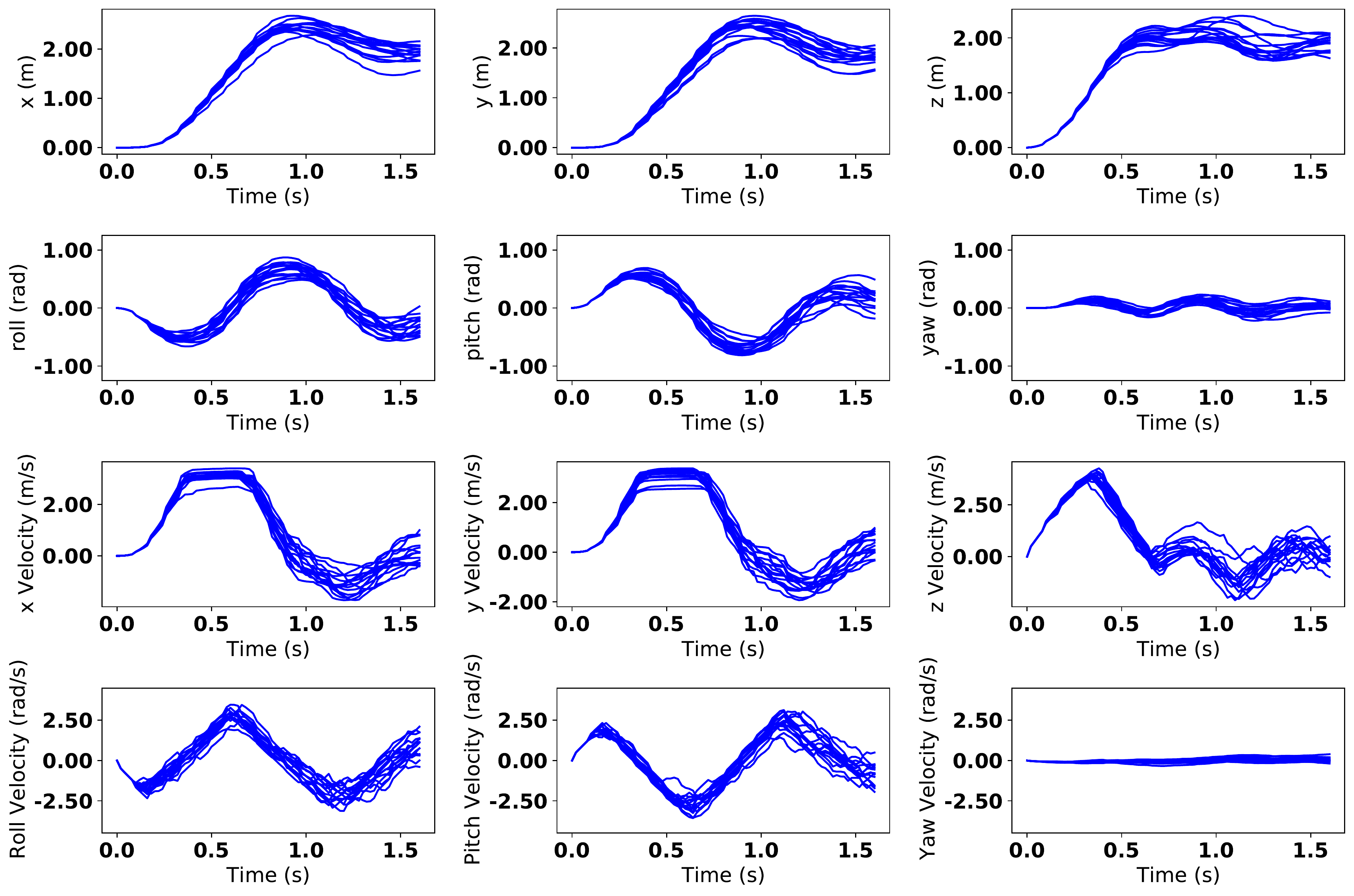}
    \caption{Trajectories from the Quadcopter problem with stochastic dynamics.
        \textit{Top}: A comparison of the posterior state from the particle filter and the ground truth from a single trajectory.
        \textit{Bottom}: Different trajectory realizations of the stochastic dynamics.}
    \label{fig:pf_control_noise/quadcopter}
\end{figure}

The results can be found in \Cref{fig:pf_quadcopter/est_state_comp} and \Cref{fig:pf_control_noise/quadcopter} respectively.

\clearpage
\newpage

\subsection{Parameters}
To alleviate the particle degeneracy problem from particle filters, artificial process noise is added to each sample. The artificial process noises added are zero mean, and their covariances are

\begin{align*}
    \Sigma_{\textrm{Pendulum\_est}} &= \textrm{diag}([\num{1e-5}, \num{1e-5}, \num{1e-9}]) \\
    \Sigma_{\textrm{Pendulum\_no\_est}} &= \textrm{diag}([0.2, 0.2, 0]) \\
    \Sigma_{\textrm{Cartpole\_est}} &= \textrm{diag}([0.001, 0.001, 0.001, 0.001, \num{1e-6}]) \\
    \Sigma_{\textrm{Cartpole\_no\_est}} &= \textrm{diag}([0.1, 0.3, 0.3, 0.2, \num{1e-6}, 0]) \\
    \Sigma_{\textrm{Quadcopter\_est}} &= \textrm{diag}([0.02, 0.02, 0.02, 0.03, 0.03, 0.03, 0.04, 0.04, 0.04, 0.04, 0.04, 0.04, 0.001]) \\
    \Sigma_{\textrm{Quadcopter\_no\_est}} &= \textrm{diag}([0.05, 0.05, 0.05, 0.03, 0.03, 0.003, 0.04, 0.04, 0.04, 0.04, 0.04, 0.04, \num{1e-9}])
\end{align*}

where the $\textrm{est}$ subscript denotes the case of uncertain parameters and $\textrm{no\_est}$ denotes the cases of uncertain initial condition and stochastic dynamics. Note that the required artificial process noise is much higher for simulation runs that do not estimate the parameters. This is because the incorrect priors for the parameters result in dynamics that are very different from ground truth, leading to the particle filter diverging very quickly.

\end{document}